\documentclass[a4paper, 11pt]{article}
\usepackage{a4}
\usepackage{amssymb}
\usepackage{amsmath}
\usepackage{epsfig}
\usepackage[english]{babel}
\usepackage{graphicx}
\usepackage{float}
\usepackage{algorithmic}
\usepackage{algorithm}
\usepackage{booktabs}

\usepackage{subfig}
\usepackage{amsmath}
\usepackage{authblk}
\usepackage{url}
\usepackage[colorlinks=true,linkcolor=blue]{hyperref}
\title{Small ensembles of kriging models for optimization}
\author[1, 2]{Hossein Mohammadi}
\author[2, 1]{Rodolphe Le Riche \thanks{Corresponding author: Ecole Nationale Supérieure des Mines de Saint Etienne, Institut H. Fayol, 158, Cours Fauriel,
42023 Saint-Etienne cedex 2 - France \\
Tel : +33477420023 \\
Email: leriche@emse.fr}}
\author[1, 2]{Eric Touboul}
\affil[1]{Ecole des Mines de Saint Etienne, H. Fayol Institute}
\affil[2]{CNRS LIMOS, UMR 5168}
\begin{document}
\date{}
\maketitle
\vskip .5cm
%=================================================================================================================================================================
\section*{Abstract}
%=================================================================================================================================================================
The Efficient Global Optimization (EGO) algorithm uses a conditional Gaussian Process (GP) to approximate an objective function known at a finite number of observation points and sequentially adds new points which maximize the Expected Improvement criterion according to the GP. The important factor that controls the efficiency of EGO is the GP covariance function (or kernel) which should be chosen according to the objective function. 
Traditionally, a parameterized family of covariance functions is considered whose parameters are learned through statistical procedures such as maximum likelihood or cross-validation. However, it may be questioned whether statistical procedures for learning covariance functions are the most efficient for optimization as they target a global agreement between the GP and the observations which is not the ultimate goal of optimization. Furthermore, statistical learning procedures are computationally expensive.
The main alternative to the statistical learning of the GP is self-adaptation, where the algorithm tunes the kernel parameters based on their contribution to objective function improvement. After questioning the possibility of self-adaptation for kriging based optimizers, this paper proposes a novel approach for tuning the length-scale of the GP in EGO: At each iteration, a small ensemble of kriging models structured by their length-scales is created. All of the models contribute to an iterate in an EGO-like fashion. Then, the set of models is densified around the model whose length-scale yielded the best iterate and further points are produced. 
Numerical experiments are provided which motivate the use of many length-scales. The tested implementation does not perform better than the classical EGO algorithm in a sequential context but show the potential of the approach for parallel implementations. 
\vskip\baselineskip
\noindent Keywords: Continuous Global Optimization; EGO; Gaussian processes; Kernel parameters; optimization based on surrogate ensembles. \\
%=================================================================================================================================================================
\section{Introduction}
%=================================================================================================================================================================
The EGO optimization algorithm uses a kriging model, which is a conditional Gaussian process (GP) \cite{GPML}, for predicting objective function values and quantifying the prediction uncertainty. 
The shapes of sample paths of a GP such as its smoothness, periodicity, etc. are controlled by the covariance function of the process, also known as its kernel. Traditionally, a parameterized family of covariance functions is considered whose parameters are estimated. 

The kernel parameters are often estimated by statistical approaches like maximum likelihood (ML)\cite{ying1991} or cross validation (CV) \cite{zhang2010}. ML and CV are compared in \cite{bachoc2013} when the covariance structure of a GP is misspecified. It is recommended in \cite{li2005} to use a penalized likelihood for the kriging models when the sample size is small. However, the efficiency of such statistical approaches, which aims at learning the objective function globally, remains questionable in the context of optimization. 
For example, in the EGO algorithm if the design points do not carry enough information about the true function, the parameters are not estimated correctly. Theses parameters are then plugged into the expected improvement (EI) criterion that may lead to disappointing results \cite{jones2001, benassi2011}. 

Not surprisingly, several methods alternative to ML and CV have been proposed to tune the kernel parameters. For instance, in \cite{frean2008} the kernel parameters are estimated with a log normal prior density assumption over them.
A fully Bayesian approach is used in \cite{benassi2011, tajbakhsh2015}. In \cite{jones1998, forrester2008}, the process of estimating parameters and searching for the optimum are combined together through a likelihood which encompasses a targeted objective. In \cite{wang2013}, the bounds on the length-scales values are changing with the iterations following an a priori schedule. 

Another drawback of statistical learning procedures such as ML and CV in the context of moderately expensive functions\footnote{We call ``moderately expensive'' functions that take between 10 seconds and an hour to be evaluated at one point.} is their computational complexity as they involve the repeated inversion of an $n \times n$ covariance matrix (where $n$ is the number of available observations) where each inversion needs of the order of $n^3$ operations.

This paper considers isotropic kernels and investigates an alternative approach to tuning the length-scale parameter. 
In this approach, a small set of length scales (hence GP models) is first tested as alternative ways to consider the objective function, independently of their statistical relevance. 
The set is completed based on the direct contribution of the best model to the optimization. 
The method is based on ensembles of surrogates. It can also be seen as weakly self-adaptive in the sense of self-adaptive algorithms \cite{Baeck-book, hansen2001} where no statistical measure intervenes in the building of the representation which the optimization algorithm has of the objective function.

Ensembles of surrogates have attracted a lot of attention from the machine learning community for prediction \cite{hess2013}, but fewer contributions seem to address surrogate ensembles for optimizing. Several approaches have been proposed that aggregate the metamodels of the ensemble into a hopefully better metamodel either by model selection or by mixing the models. This better metamodel is then used by the optimization algorithm \cite{acar2009, Chaudhuri_ensemble2013, goel2007}.

On the opposite, other previous optimization methods take advantage of all the metamodels in the set as a diversity preserving mechanism (in addition to, of course, a way to reduce the number of calls to the objective function), in the context of evolutionary computation \cite{jin2004, lu2013} or more generally \cite{viana2013}. The algorithm studied in this text belongs to this category.

Another classification can be made with respect to the homogeneity (all metamodels are of the same type) or heterogeneity of the ensemble. There has been recent contributions to optimization algorithms that rely on a homogeneous set of kriging models: in \cite{kleijnen2014} the ensembles are built by bootstrap on the data and serve as a way to estimate model uncertainty for later use in optimization; in \cite{villemonteix2008}, the metamodels are the trajectories of a GP and their contributions are aggregated through an uncertainty reduction criterion (on the entropy of the global optima of the trajectories).  
The optimization algorithm investigated here also relies on an homogeneous ensemble of GP models.
%=================================================================================================================================================================
\section{EGO algorithm overview}
\begin{figure}[htpb] 
\centering
\includegraphics[width=0.49\textwidth]{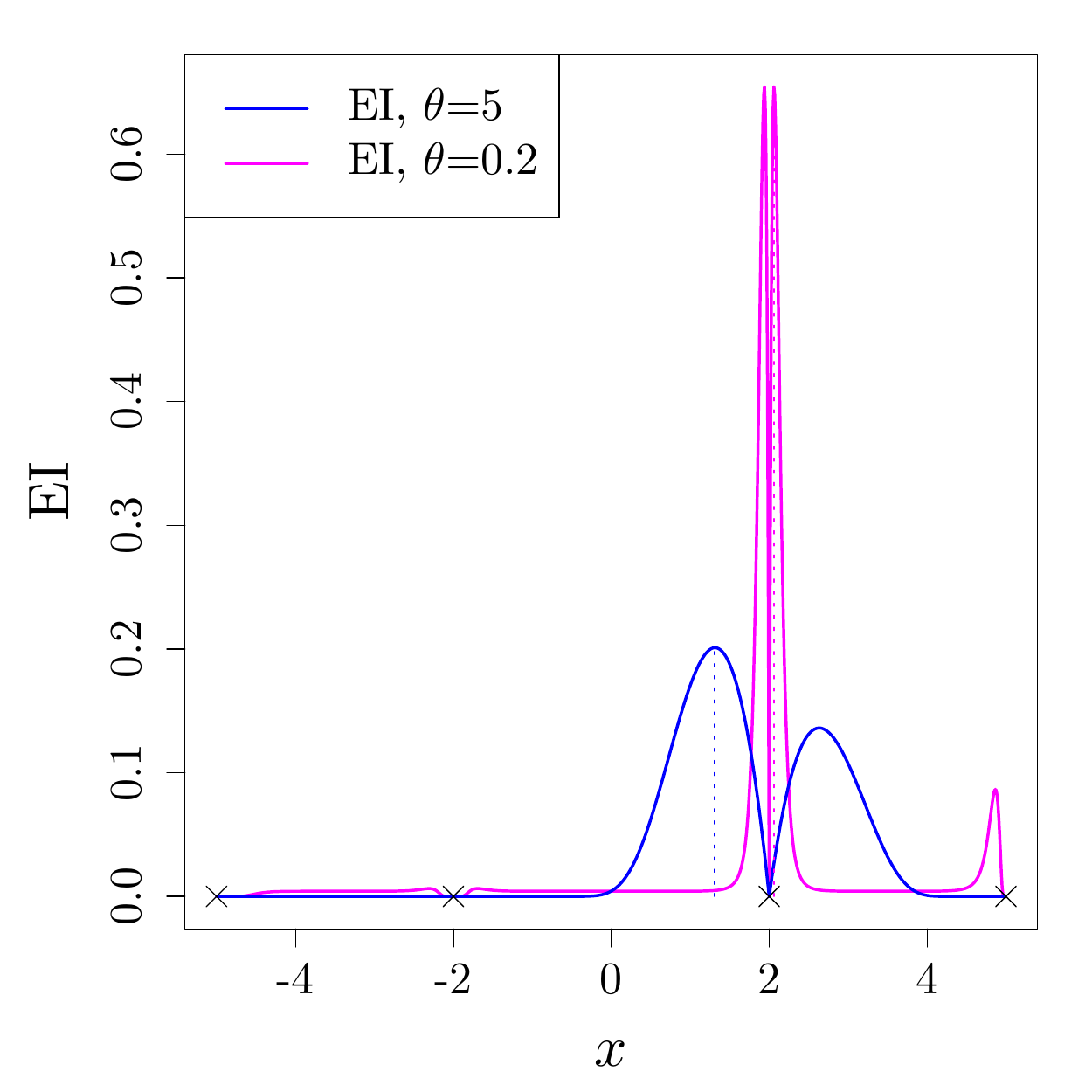}
\includegraphics[width=0.49\textwidth]{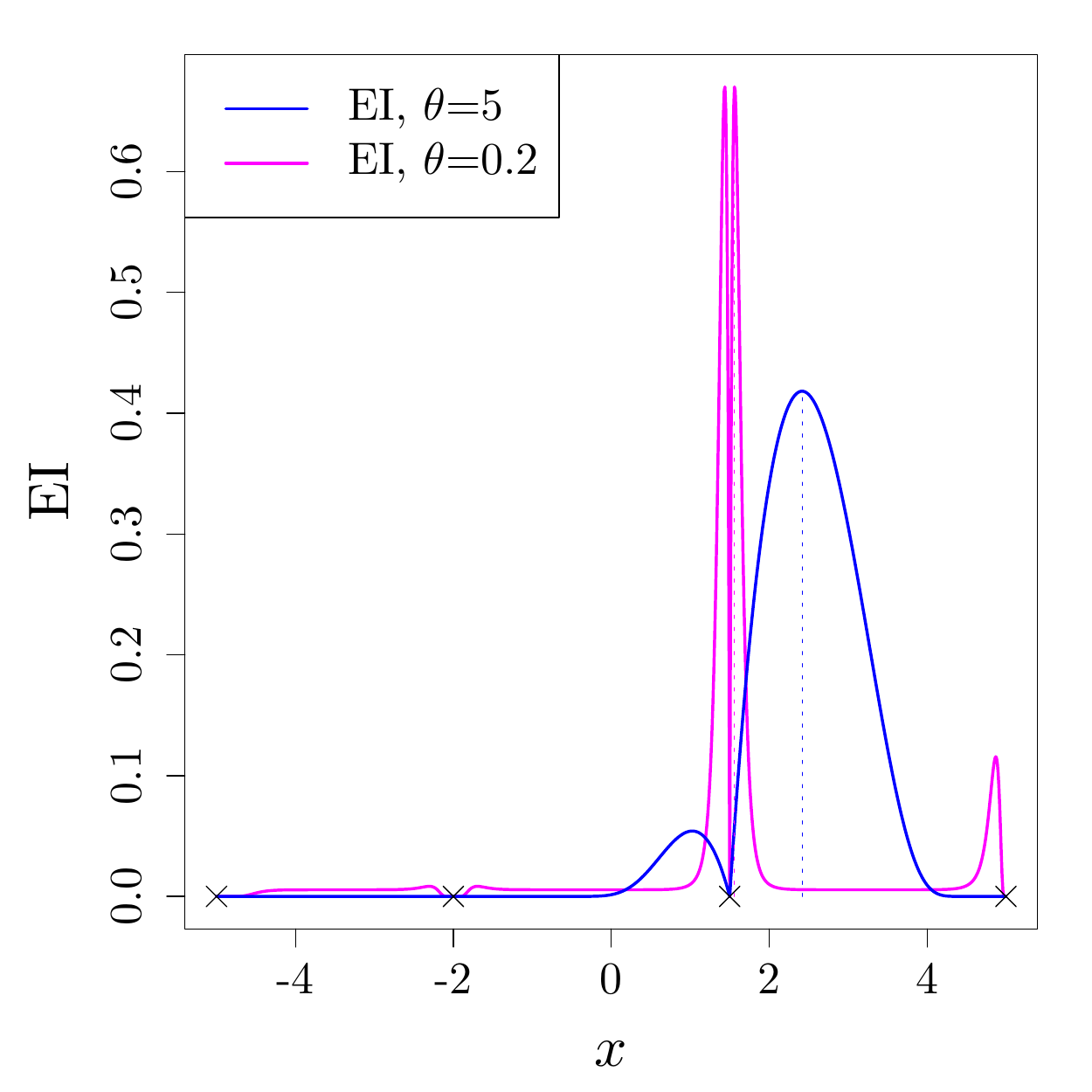}
\caption{Effect of DoE and length-scale on EI function. The function to be optimized is the Sphere whose global minimum is located at $2.5$. The blue and magenta curves represent the EI of kriging models with length-scales equal to $5$ and $0.2$, respectively. The crosses indicate the location of design points. The other parameters are fixed. The location of the third sample point changes from $2$ to $1.5$ in the right picture. 
} 
\label{theta_effect_on_EI}
\end{figure}

%=================================================================================================================================================================
EGO is a sequential model-based optimization algorithm. It starts with an initial design of experiments (DoE). At each iteration, one point which maximizes the Expected Improvement (EI) according to the current kriging model is added to the DoE. Then, the kernel parameters are re-estimated and the kriging model is updated.
\begin{algorithm}
\caption{Efficient Global Optimization Algorithm (EGO)}
\label{EGO_algorithm}
\begin{algorithmic}
\STATE Create an initial design: $\textbf{X} = \left[ \textbf{x}^1, \dots , \textbf{x}^n \right]^T$.
\STATE Evaluate the functions at $\textbf{X}$, $\textbf{y} = f(\textbf{X})$.
\STATE Fit a kriging model to the data points $(\textbf{X}, \textbf{y})$ = estimate $\theta, \mu, \sigma^2$
\WHILE{\NOT stop}
\STATE $\textbf{x}^{n+1} ~\leftarrow ~ \arg\max_{\textbf{x} \in \mathcal S}~EI(\textbf{x})$ and add $\textbf{x}^{n+1}$ to $\textbf{X}$.
\STATE $y^{n+1} ~\leftarrow ~ f(\textbf{x}^{n+1})$ and add $y^{n+1}$ to $\textbf{y}$.
\STATE Re-estimate the parameters ($\theta, \mu, \sigma^2$) and update the kriging model.
\ENDWHILE
\end{algorithmic}
\end{algorithm}
%=================================================================================================================================================================
The location of $\textbf{x}^{n+1}$, where $\textbf{x}^{n+1} = \arg \max_{x \in \mathcal S} EI(\textbf{x})$,  depends on the current DoE, $\textbf{X}$, $\textbf{y}$, the kriging trend, $\mu$, and the kernel parameters: the length-scale, $\theta$, and the process variance, $\sigma^2$. We use $\textbf{x}^{n+1} = g(\textbf{X}, \mu, \theta, \sigma^2)$ to denote that $\textbf{x}^{n+1}$ is a function of the above-mentioned parameters. 
Figure \ref{theta_effect_on_EI} illustrates how the DoE and the magnitude of length-scale affect the EI. 
%%%%%%%%%%%%%%%%%%%%%%%%%%%%%%%%%%%%%%%%%%%%%%%%%%%%%%%%%%%%%%%%%%%%%%%%%%%%%%%%%%%%%%%%%%%%%%
%%%%%%%%%%%%%%%%%%%%%%%%%%%%%%%%%%%%%%%%%%%%%%%%%%%%%%%%%%%%%%%%%%%%%%%%%%%%%%%%%%%%%%%%%%%%%%
 
\begin{figure}[H] 
\centering
\includegraphics[width=0.49\textwidth]{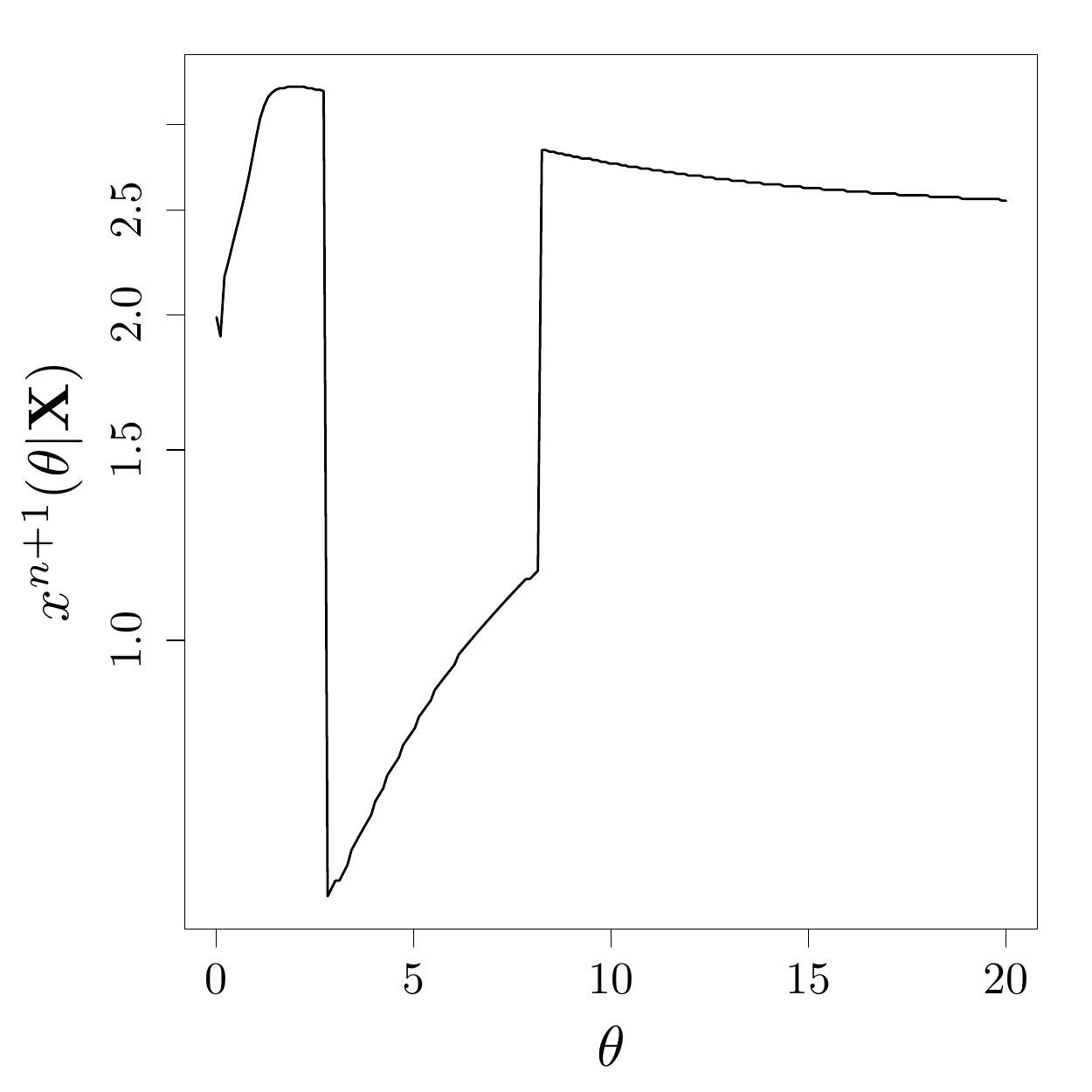}
\includegraphics[width=0.49\textwidth]{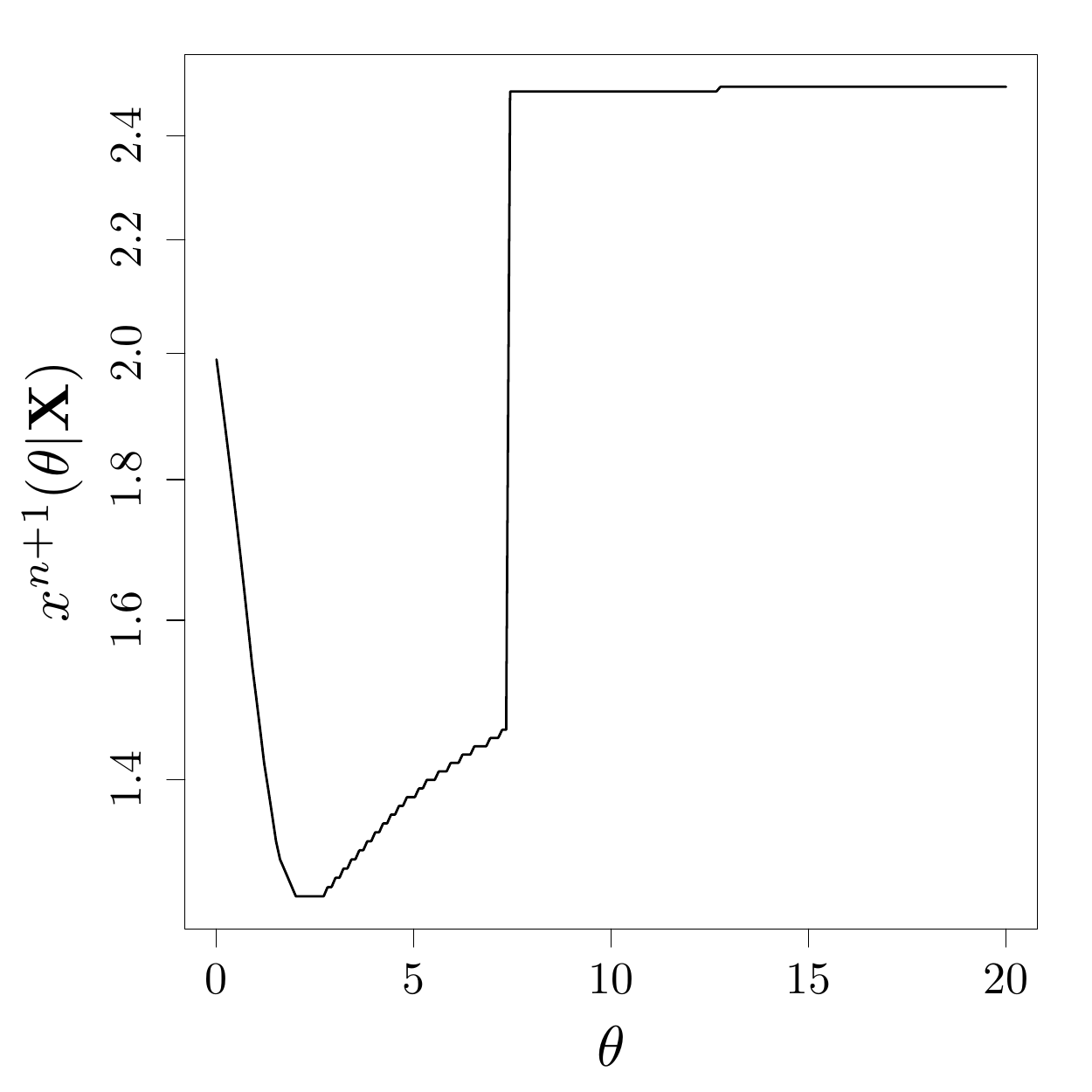}
\caption{Illustration of all possible next infill sample points with $\textbf{X} = \{-5, -2, 2, 5 \}$ as the DoE. The true functions are Sphere (left, as in Figure~\ref{theta_effect_on_EI}) and Ackley (right) in dimension $1$. For $\theta$ values larger than, say $\theta \geq 8$, the location of $x^{n+1}$ is quite stable and close to 2.5, the location of the global minimum. While large $\theta$'s lead to the global optimum of the Sphere for any $\textbf X$, it is a coincidence for Ackley's function.} 
\label{NextPoint_as_fun_theta}
\end{figure}

Among the parameters of the EI criterion, $\textbf{X}$ and $\theta$ play a prominent role 
because once $\textbf{X}$ and $\theta$ are fixed, the ML estimations of $\mu$ and $\sigma^2$ have a closed-form expression \cite{GPML}:  
\begin{eqnarray}
\hat{\mu} =& \frac{\textbf{1}^\top \textbf{R}^{-1}(\theta)\textbf{y}}{\textbf{1}^\top \textbf{R}^{-1}(\theta)\textbf{1}}, \\
\label{kriging_trend_estim}
\hat{\sigma}^2 =& \frac{(\textbf{y} - \hat{\mu} \textbf{1})^\top \textbf{R}^{-1} (\theta) (\textbf{y} - \hat{\mu} \textbf{1})}{n}.
\label{process_variance_estim}
\end{eqnarray}
Accordingly, $\textbf{x}^{n+1}$ can be expressed as a function of $\textbf{X}$ and $\theta$. For example, Figure \ref{NextPoint_as_fun_theta} shows all plausible next infill sample points by changing the length-scale for a given DoE. 
%%%%%%%%%%%%%%%%%%%%%%%%%%%%%%%%%%%%%%%%%%%%%%%%%%%%%%%%%%%%%%%%%%%%%%%%%%%%%%%%%%%%%%%%%%%%%%
%%%%%%%%%%%%%%%%%%%%%%%%%%%%%%%%%%%%%%%%%%%%%%%%%%%%%%%%%%%%%%%%%%%%%%%%%%%%%%%%%%%%%%%%%%%%%%
%=================================================================================================================================================================
\section{Tuning the length-scale from an optimization point of view: a study on self-adaptation}
%=================================================================================================================================================================
When the kernel parameters are estimated by ML, the selected kriging model has statistical best agreement with the observed data. However, the goal of using EGO, like other optimization algorithms, is to solve an optimization problem with the least number of function evaluations. 
In other words, the main goal is the fast convergence of EGO even if the kriging model does not represents well the true function.
This idea is similar to the notion of ``self-adaptation'' in evolutionary optimization \cite{Baeck-book,hansen2001}.

To investigate the potential of tuning the length-scale $\theta$ in an optimization oriented, greedy, self-adaptive way, we first tested a 
theoretical algorithm that tries all possible values of $\theta$ in the range $[0.01,20]$. The true objective functions of the points 
that maximize the expected improvement for each of these length-scale $\theta$ value is calculated, 
$\textbf{x}^{n+1}(\theta|\textbf{X})=\arg \max_{x \in \mathcal S} EI(\textbf{x};\theta)$.  
This makes this algorithm not practical in the context of 
expensive problems. Then, the iterate associated to the best objective function, \newline $\textbf{x}^{\text{sel}} = \arg \min_{\textbf{x}^{n+1}} f(\textbf{x}^{n+1}(\theta|\textbf{X}))$, is added to the Design of Experiment $\textbf X$, the kriging model is updated, and the algorithm loops. This algorithm is sketched in the flow chart \ref{EGO_greedy}.

\begin{figure}[H] 
\centering
\includegraphics[width=0.49\textwidth]{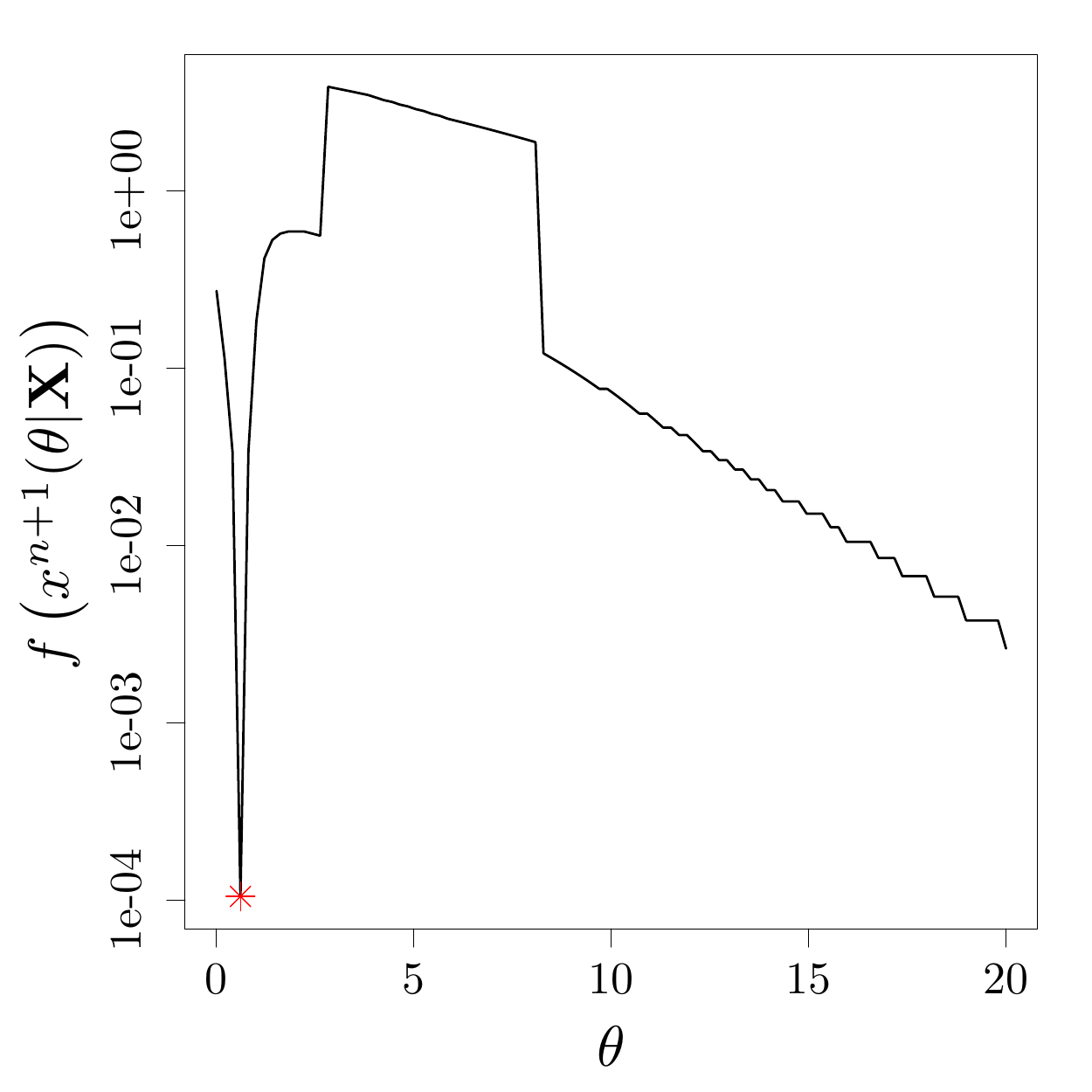}
\includegraphics[width=0.49\textwidth]{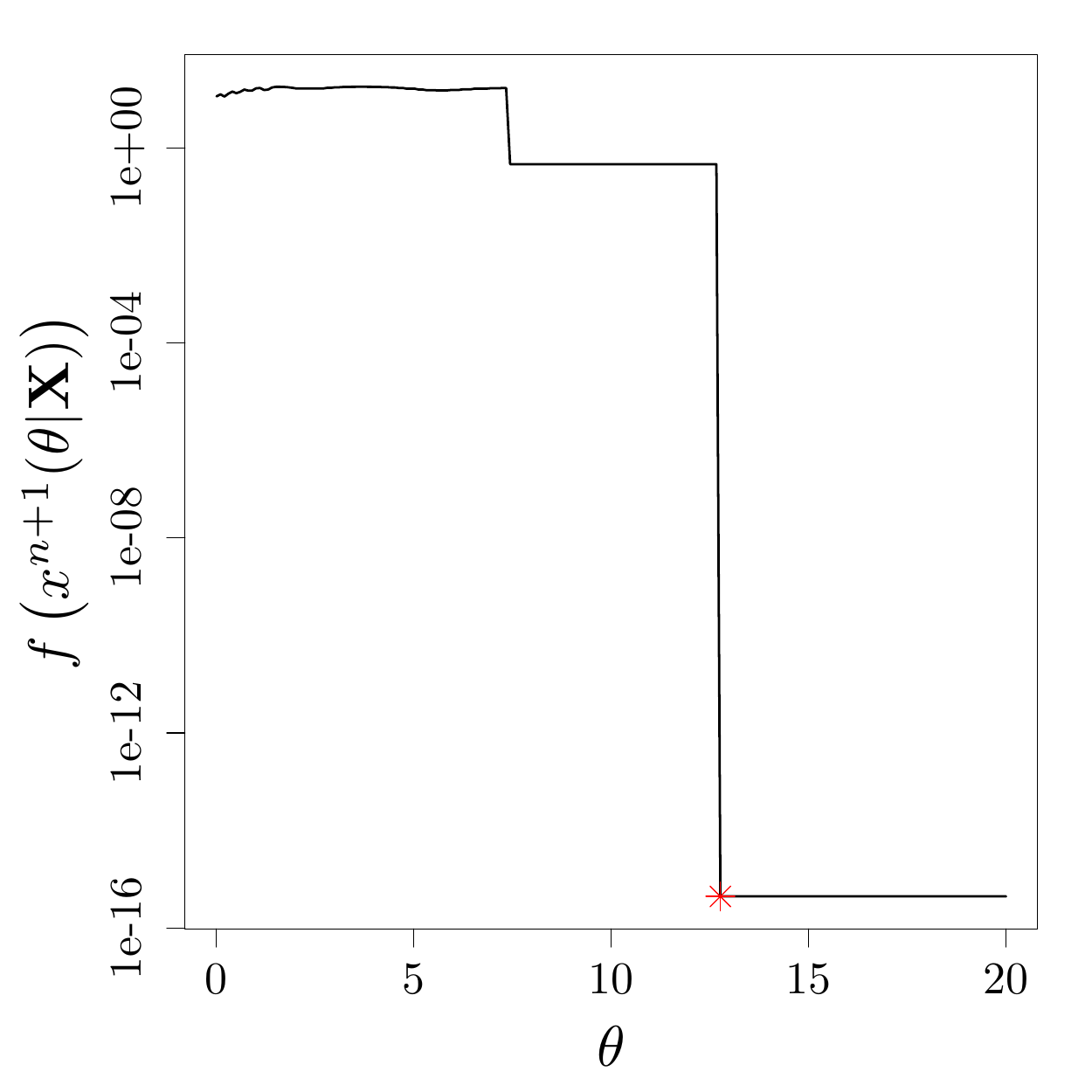}
\caption{Function values of $x^{n+1}$ already shown in Figure \ref{NextPoint_as_fun_theta}. The asterisk indicate the correlation length-scale, $\theta^*$, which causes the maximum improvement in the objective function. In this example, $\theta^*$ is different from $\hat{\theta}_{ML}$, estimated by ML,: $\theta^* = 0.61271$ and $\hat{\theta}_{ML} = 5.34$ (Sphere; left), $\theta^* = 12.7674$ and $\hat{\theta}_{ML} = 0.01$ (Ackley; right).   
Both functions have their global minimum at 2.5 and the DoE is $\textbf{X} = \{-5, -2, 2, 5 \}$.} 
\label{TrueFunc_as_fun_theta}
\end{figure}

From a one step ahead optimization point of view, the ``best'' length-scale, denoted by $\theta^*$, is the one that yields the next infill sample with the lowest objective function value, $\theta^* = \arg \min_\theta f(\textbf{x}^{n+1}(\theta|\textbf{X})) $. In the examples provided in Figure \ref{TrueFunc_as_fun_theta}, the best length-scales are shown for the two test functions (Ackley and Sphere). In this example, the best length-scales are different from the length-scales estimated by ML, see the caption of Figure \ref{TrueFunc_as_fun_theta}. 

\begin{algorithm}[H]
\caption{Toy EGO with greedy $\theta$ tuning}
\label{EGO_greedy}
\begin{algorithmic}
\STATE Create an initial design: $\textbf{X} = \left[ \textbf{x}^1, \dots , \textbf{x}^n \right]^T$
\STATE Evaluate the functions at $\textbf{X}$, $\textbf{y} = f(\textbf{X})$
\STATE Select $\textbf{x}^j~, 1 \leq j \leq n,$ for which $f(\textbf{x}^j)=\max(\textbf{y})$
\WHILE{\NOT stop (typically a limit on budget)}  
\STATE Set $\textbf{x}^{\text{sel}} = \textbf{x}^j$
\FOR{$\theta_i \in [\theta_{min}, \dots , \theta_{max}]$}
\STATE $\textbf{x}^{n+1}(\theta_i|\textbf{X}) = \arg\max_{x \in \mathcal S} EI(\textbf{x}; \theta_i)$
\IF {$f\left(\textbf{x}^{n+1}(\theta_i|\textbf{X})\right) < f(\textbf{x}^{\text{sel}})$}
\STATE $\textbf{x}^{\text{sel}}~\leftarrow ~\textbf{x}^{n+1}(\theta_i|\textbf{X})$
\ENDIF
\ENDFOR
\STATE $\textbf{X} ~\leftarrow ~ \textbf{X} \cup \textbf{x}^{\text{sel}}$
\ENDWHILE
\end{algorithmic}
\end{algorithm}

%%%%%%%%%%%%%%%%%%%%%%%%%%%%%%%%%%%%%%%%%%%%%%%%%%%%%%%%%%%%%%%%%%%%%%%%%%%%%%%%%%%%%%%%%%%%%%
%%%%%%%%%%%%%%%%%%%%%%%%%%%%%%%%%%%%%%%%%%%%%%%%%%%%%%%%%%%%%%%%%%%%%%%%%%%%%%%%%%%%%%%%%%%%%%

\begin{figure}[H] 
\centering
\includegraphics[width=0.49\textwidth]{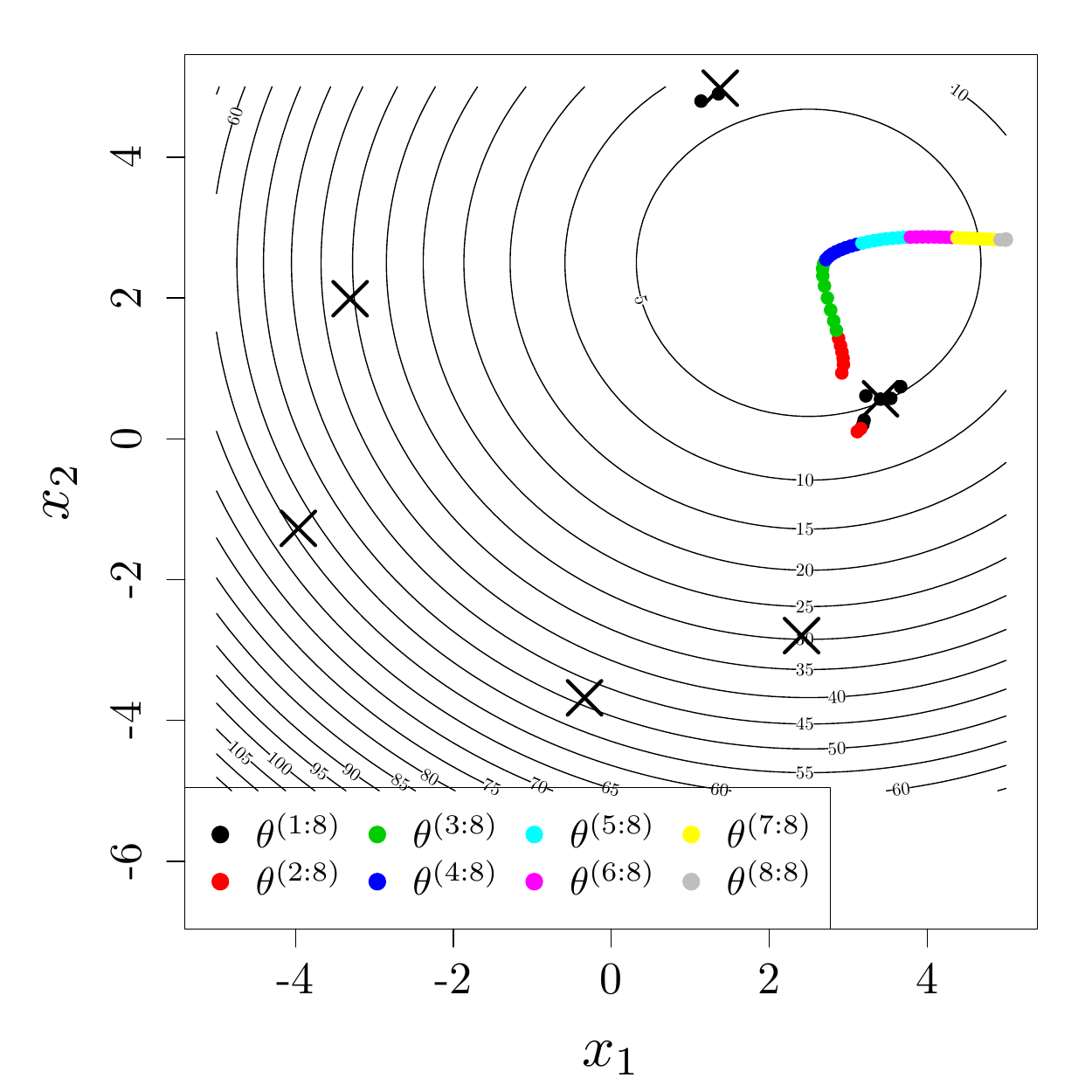}
\includegraphics[width=0.49\textwidth]{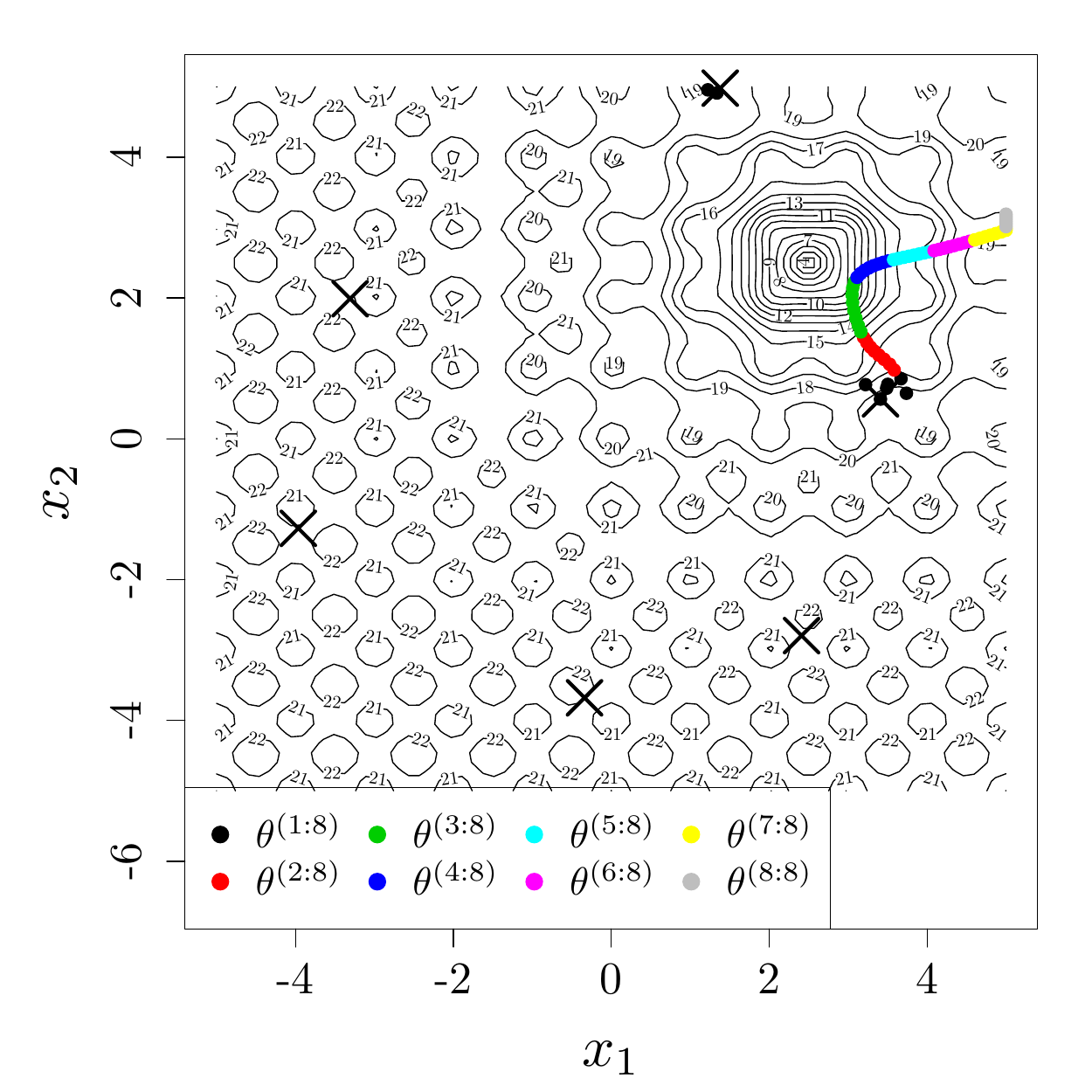}
\includegraphics[width=0.49\textwidth]{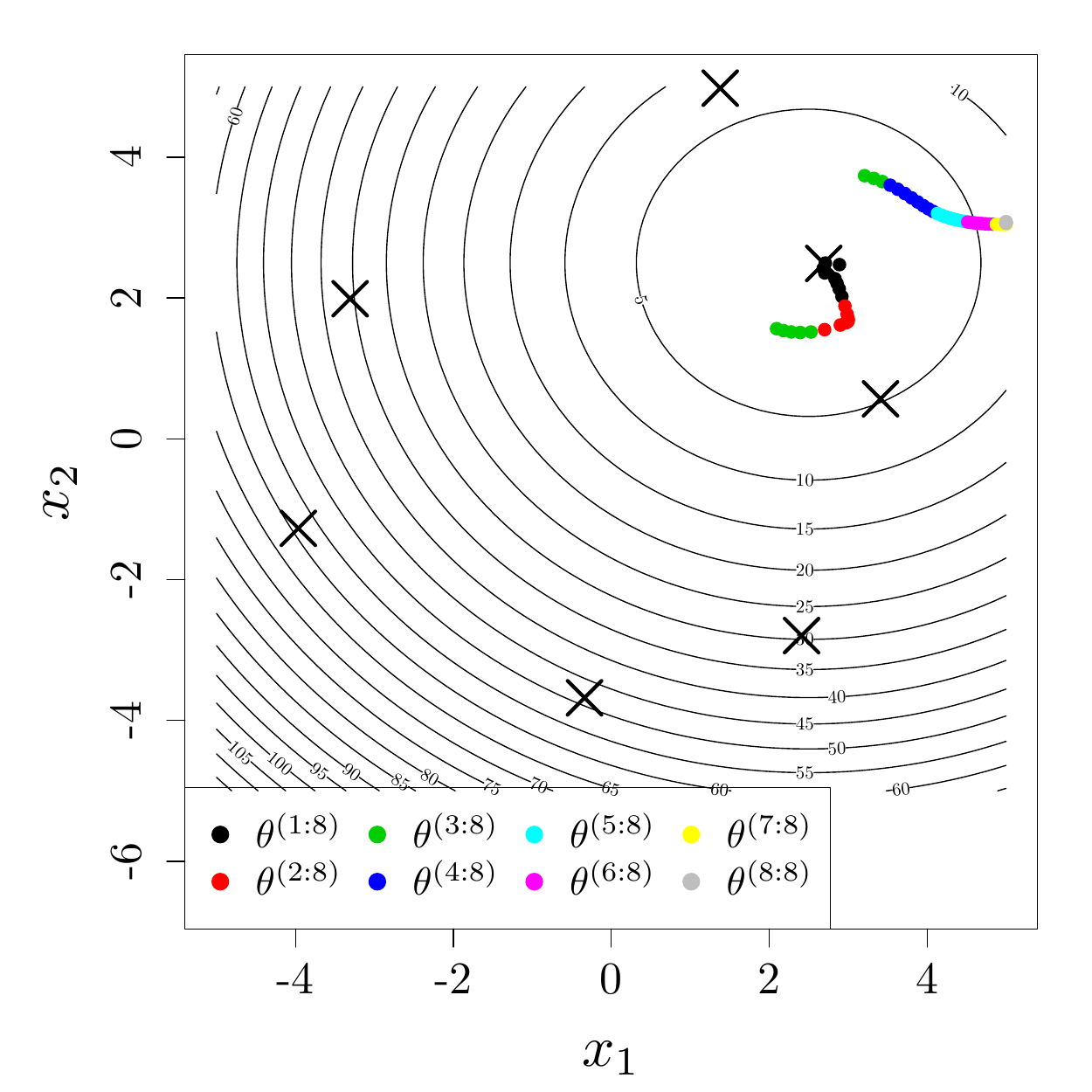}
\includegraphics[width=0.49\textwidth]{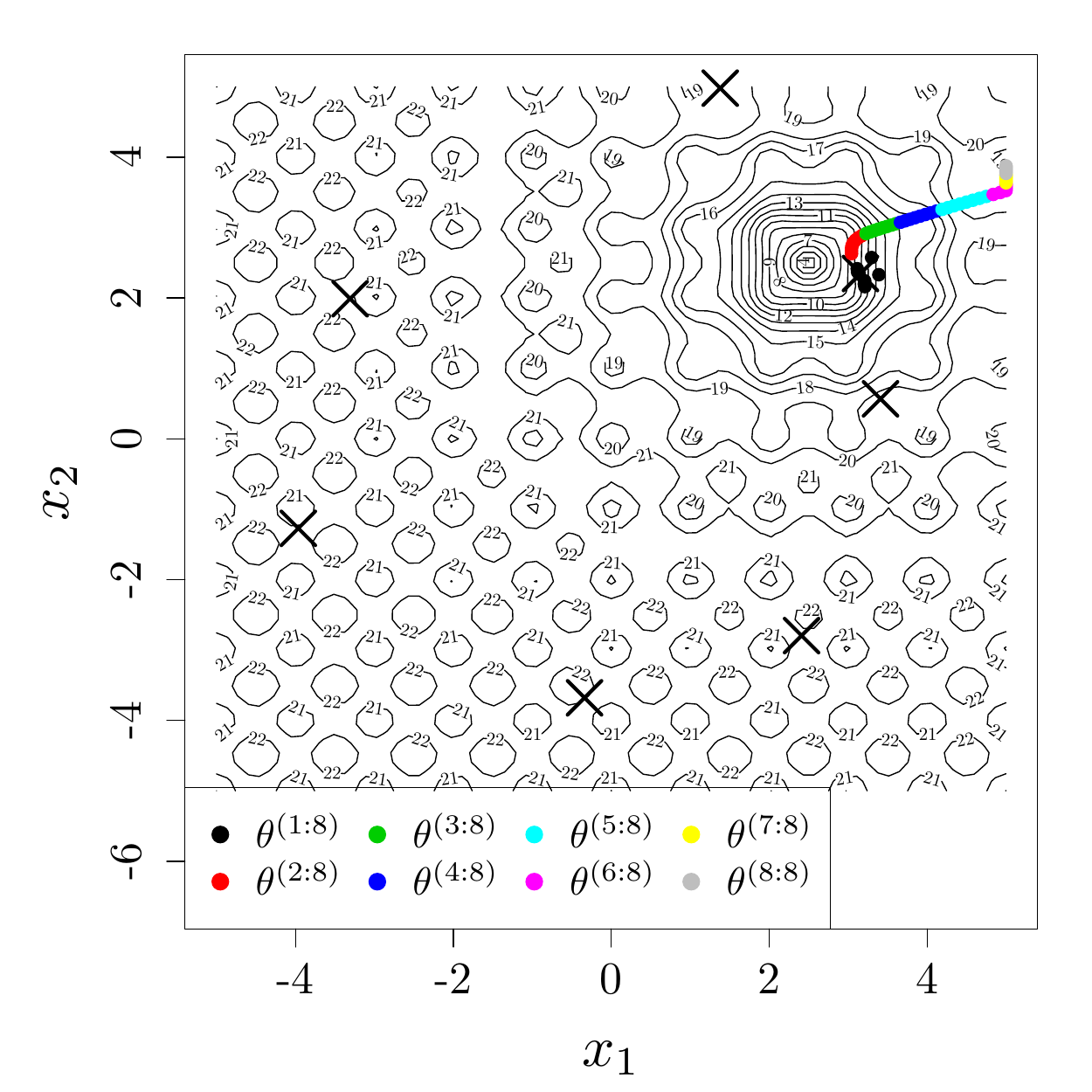}
\caption{First (top row) and second (bottom row) iterations of EGO in which $\textbf{x}^{n+1}(\theta^*|\textbf{X})=\arg \max_{x \in \mathcal S} EI(\textbf{x}|\theta^*)$ is added to the existing DoE, the crosses, on the Sphere (left) and the Ackley (right) functions. 64 equally distant length-scales are grouped into eight equal sized intervals, $\theta^{(i:8)}, i=1, \dots, 8$. The infill sample points obtained by the length-scales of a
particular group have identical color.} 
\label{x_theta_2D}
\end{figure}

We now analyze this approach in more details by providing some examples in $2D$. Figure \ref{x_theta_2D} illustrates the first and the second iterations of this algorithm again on the Sphere and Ackley functions. 
In this Figure, the location of the points that maximize the expected improvement for different length-scale values is plotted on the top of the true function contour lines. 
In total, 64 length-scales, started from 0.01, are used. The length-scales are divided into eight groups. Each group consists of eight length-scales in ascending order. The $i$th group is denoted by $\theta^{(i:8)}, i=1, \dots, 8$ and is defined as $\left[0.01 + 8(i-1)\times \alpha_{\text{increment}}, 0.01 + 8i \times \alpha_{\text{increment}} \right)$ where $\alpha_{\text{increment}} \approx  0.1$.
The infill sample points obtained by the length-scales of a particular group have identical color, see the legend of Figure \ref{x_theta_2D}.    

The first remark that can be done, and which motivates this study, is that the points visited as $\theta$ changes make a one dimensional manifold (obviously since it is parameterized by the scalar $\theta$), continuous by parts and, most interestingly, often curved towards the global optimum of the function. The discontinuities of the trajectory are associated to changes of basin of attraction during the maximization of the expected improvement.
This simple observation, even though only based on a few cases, is a hint that the volume search of global optimization algorithms might be iteratively transformed into a one dimensional search in $\theta$, with potentials for containing the ``curse of dimensionality'' (the geometric increase in search space volume as the number of dimensions increases). The difficulties of the associated problem and a possible implementation will be discussed in the next section.

\begin{figure}[H] 
\centering
\includegraphics[width=0.32\textwidth]{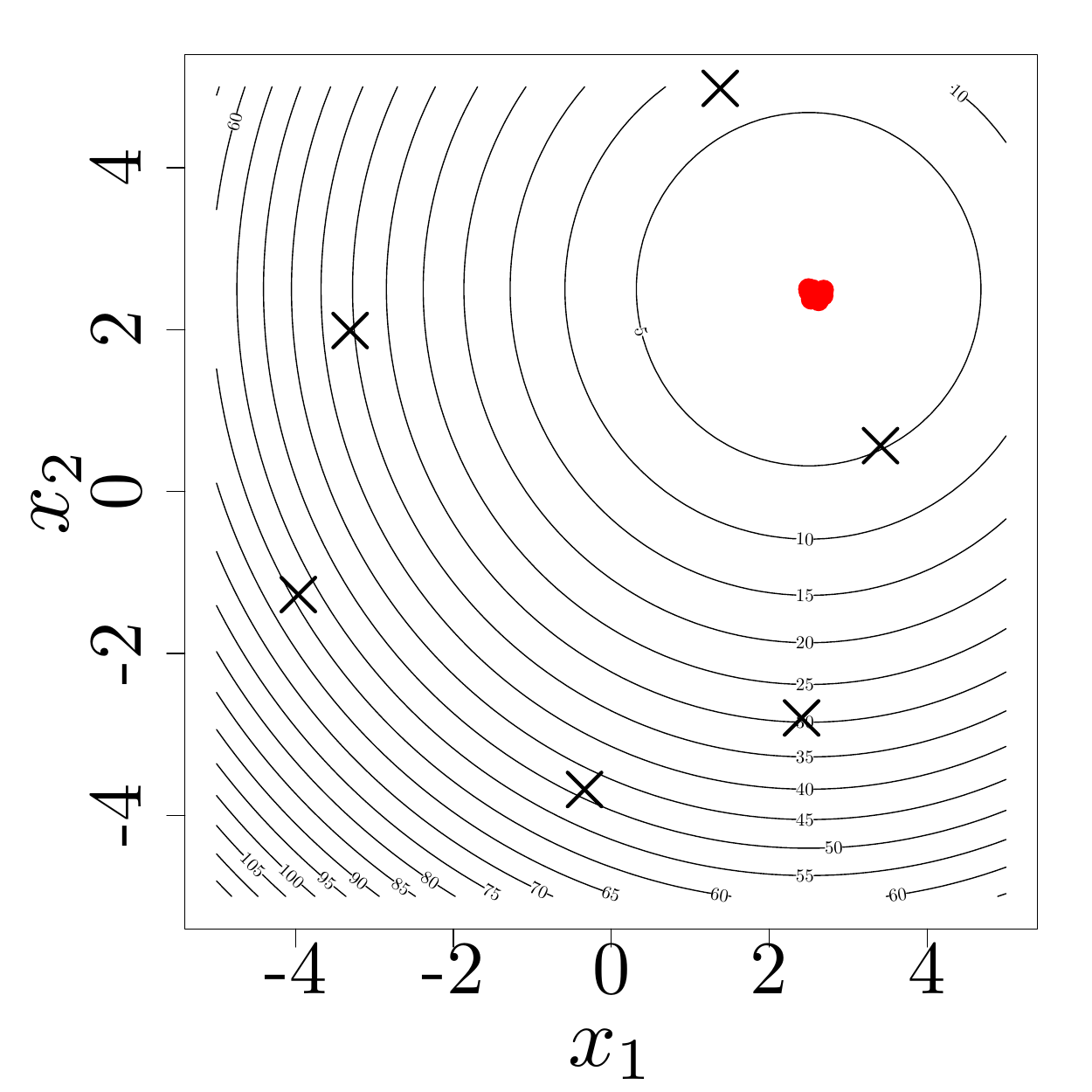}
\includegraphics[width=0.32\textwidth]{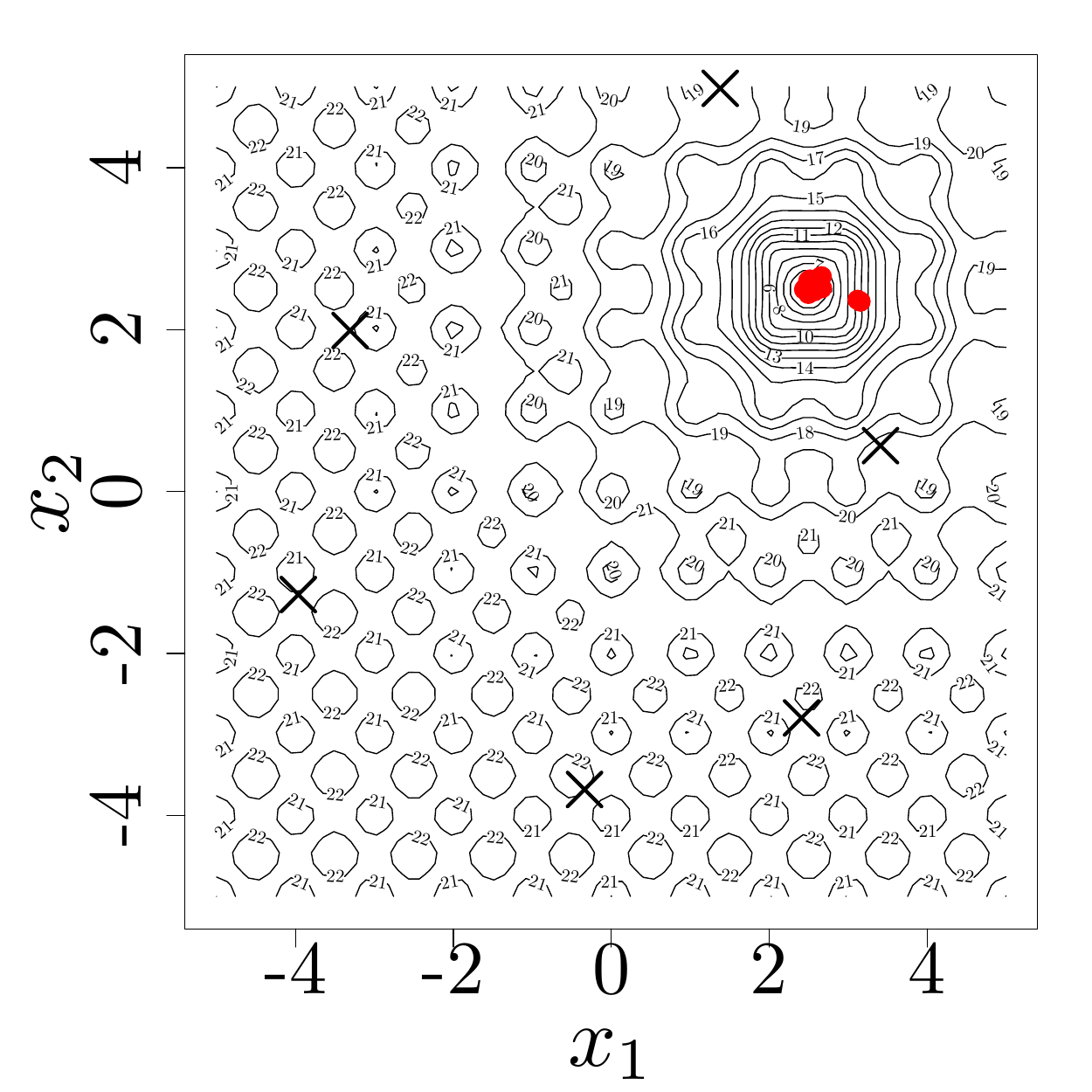}
\includegraphics[width=0.32\textwidth]{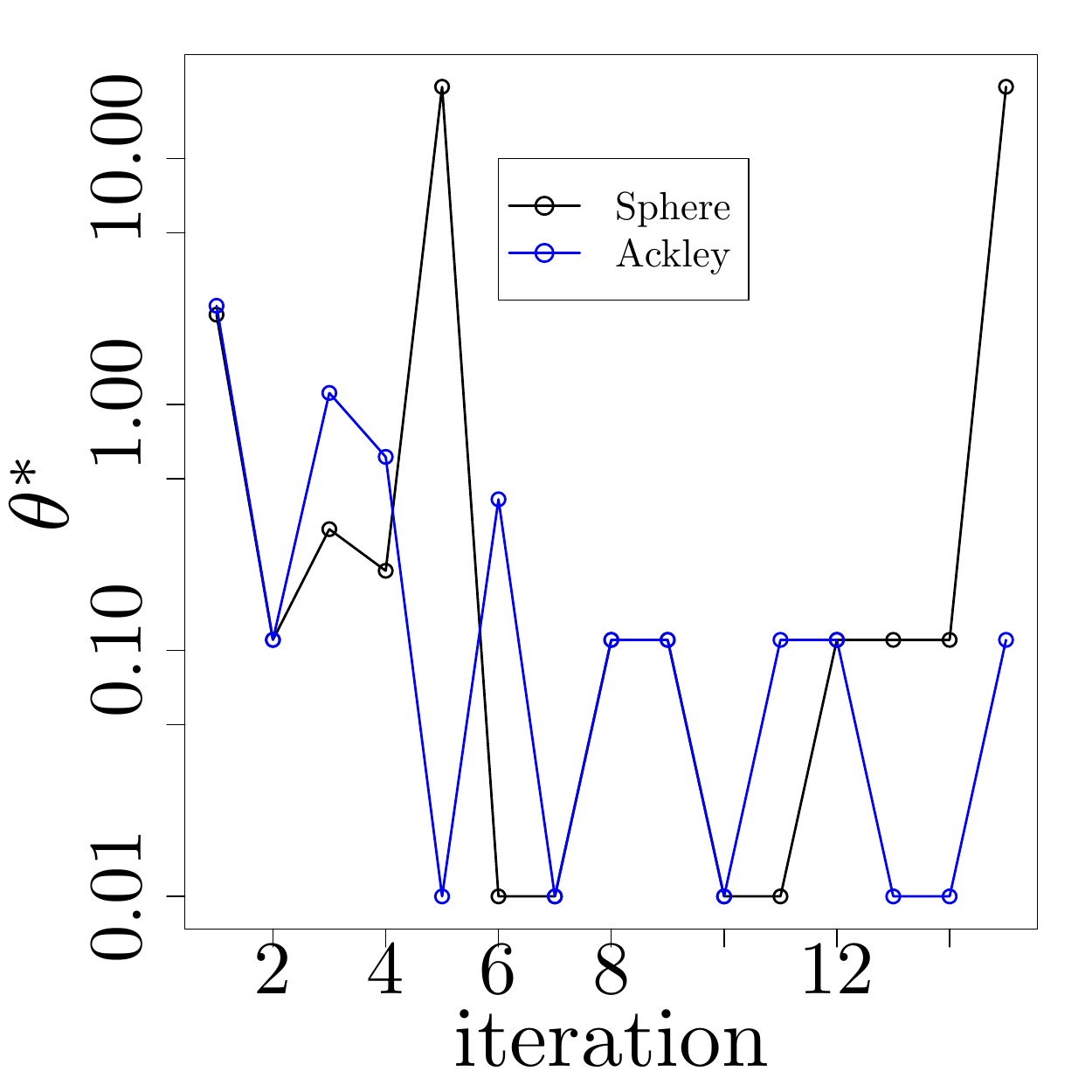}
\caption{DoEs created by the toy greedy algorithm \ref{EGO_greedy} after 15 iterations on the Sphere (left) and the Ackley (middle) functions. Right: plot of ``best'' length-scale, $\theta^{*}$. $\theta^{*}$ oscillates during optimization iterations and usually has a small magnitude after the first iterations. The y-axis is in logarithmic scale.} 
\label{15_iter}
\end{figure}

In Figure \ref{x_theta_2D}, it can be seen that the magnitude of the ``best'' length-scale in the first iteration is between 2 and 3, i.e., $\theta^* \in \theta^{(3:8)} ~~ or ~~ \theta^{(4:8)}$.
While EGO with a small length-scale samples near the best observed point (cf. the black points), EGO with large length-scale is more explorative (see yellow and grey points) \cite{mohammadi2015}. 
The search points and the length-scales obtained by the algorithm after 15 iterations are given in Figure \ref{15_iter}. It can be observed that, after the first iterations where the ``best'' length-scale magnitude, $\theta^{*}$, is of order 1, $\theta^{*}$ oscillates at usually small values. Because $\theta^{*}$ oscillates, self-adaptive strategies and Bayesian strategies based on assuming a prior density over the length-scale may not be a good strategy for optimization (at least if $\theta^*$ makes an efficient strategy). 
%%%%%%%%%%%%%%%%%%%%%%%%%%%%%%%%%%%%%%%%%%%%%%%%%%%%%%%%%%%%%%%%%%%%%%%%%%%%%%%%%%%%%%%%%%%%%%
%%%%%%%%%%%%%%%%%%%%%%%%%%%%%%%%%%%%%%%%%%%%%%%%%%%%%%%%%%%%%%%%%%%%%%%%%%%%%%%%%%%%%%%%%%%%%%

%%%%%%%%%%%%%%%%%%%%%%%%%%%%%%%%%%%%%%%%%%%%%%%%%%%%%%%%%%%%%%%%%%%%%%%%%%%%%%%%%%%%%%%%%%%%%%
%%%%%%%%%%%%%%%%%%%%%%%%%%%%%%%%%%%%%%%%%%%%%%%%%%%%%%%%%%%%%%%%%%%%%%%%%%%%%%%%%%%%%%%%%%%%%%

\begin{figure}[H] 
\centering
\includegraphics[width=0.49\textwidth]{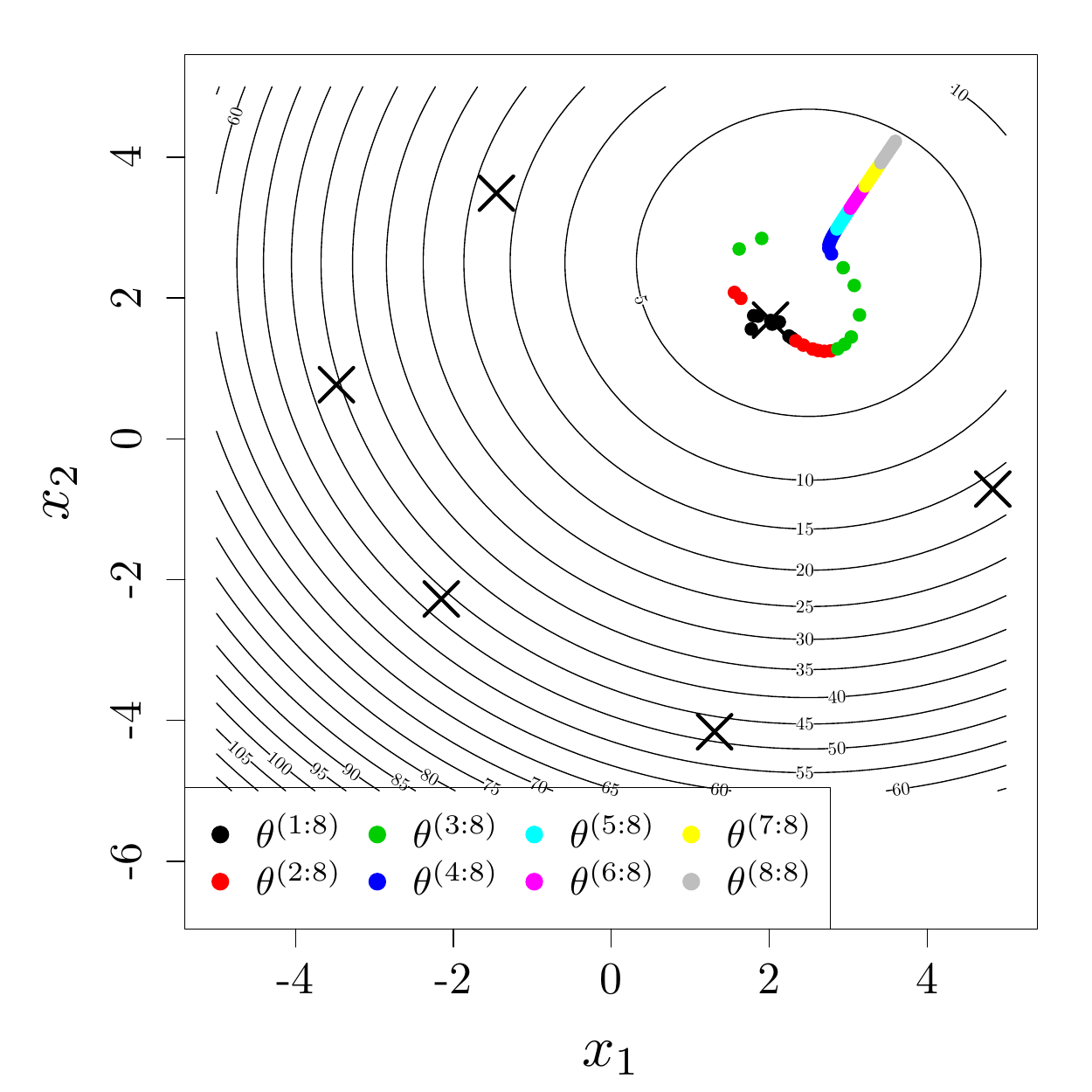}
\includegraphics[width=0.49\textwidth]{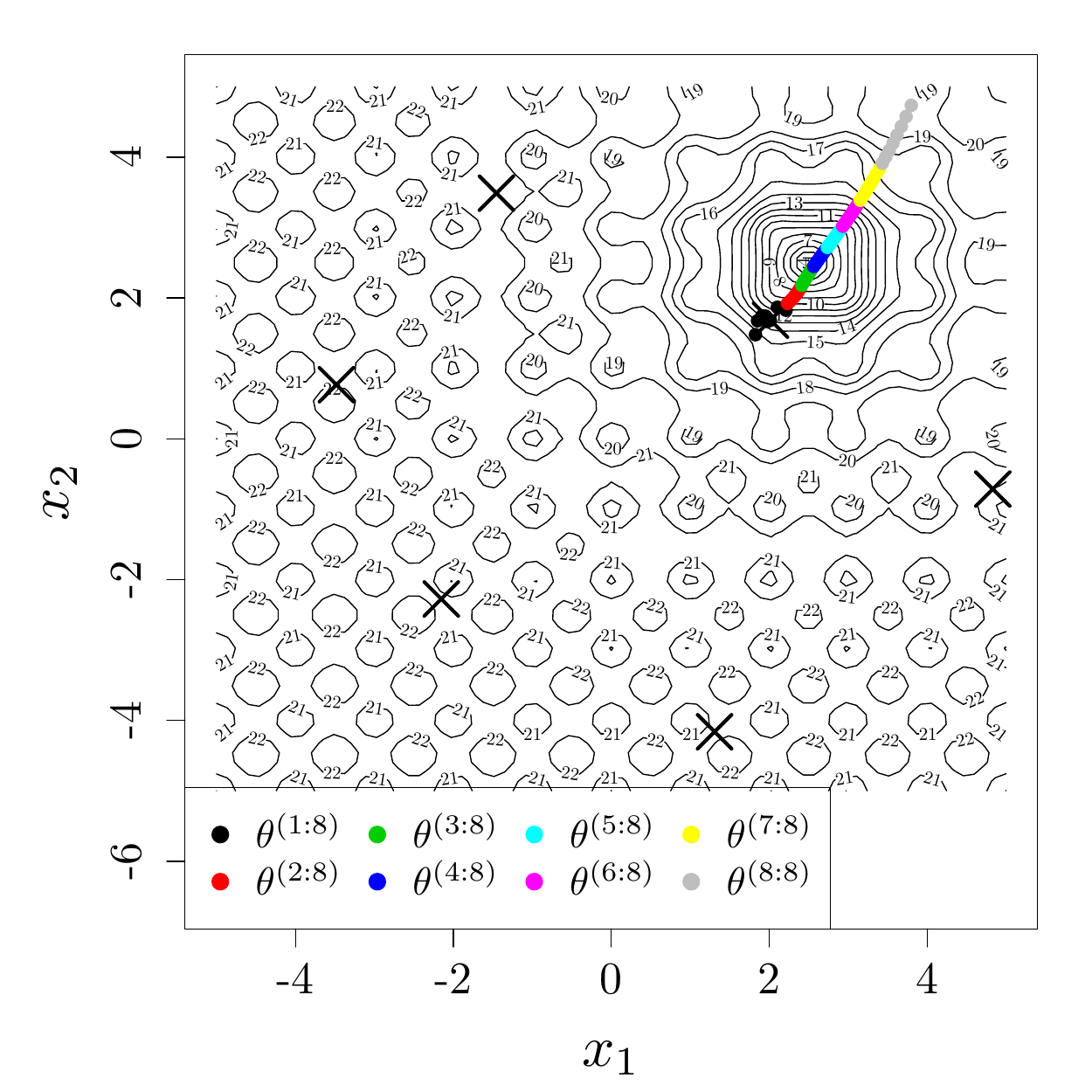}
\includegraphics[width=0.49\textwidth]{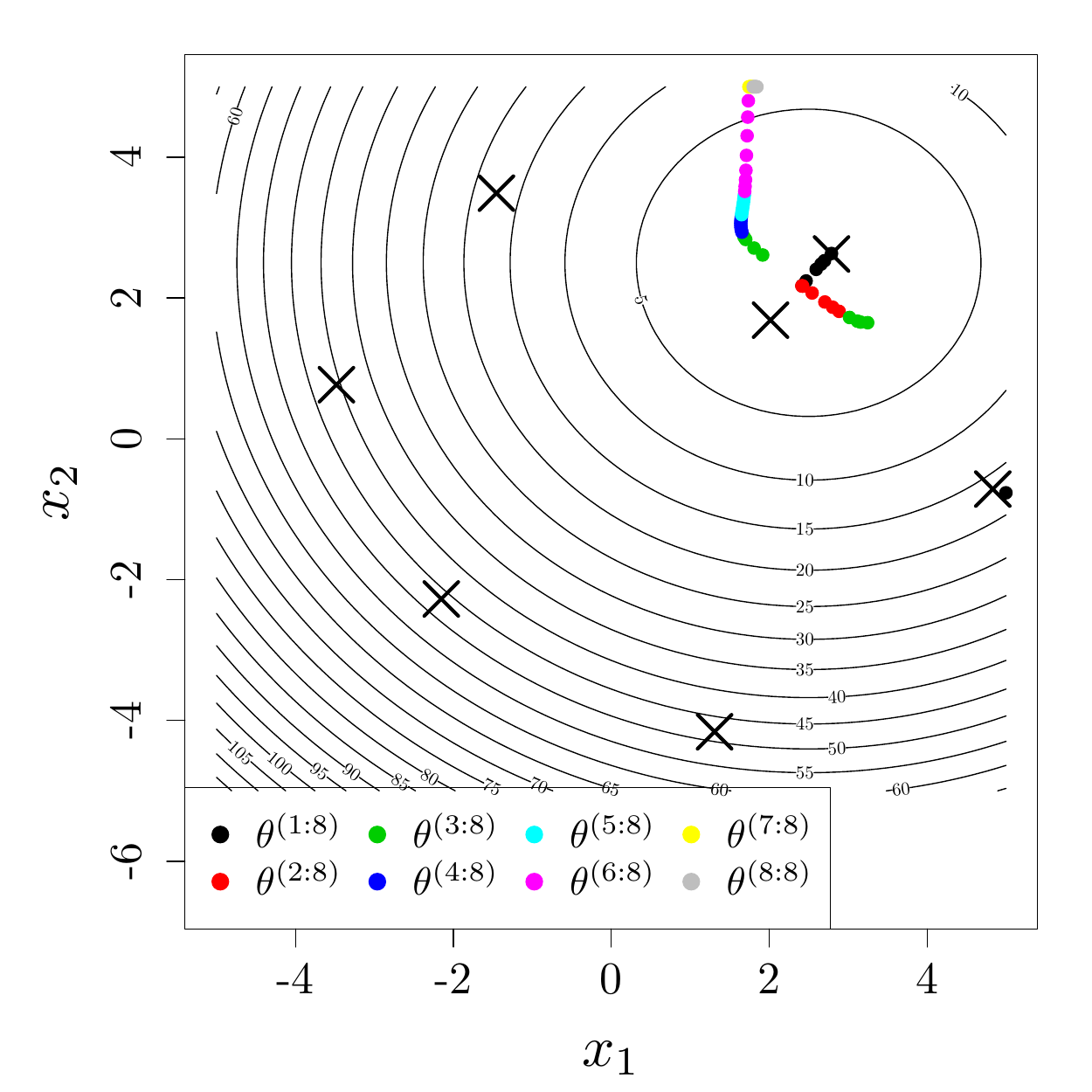}
\includegraphics[width=0.49\textwidth]{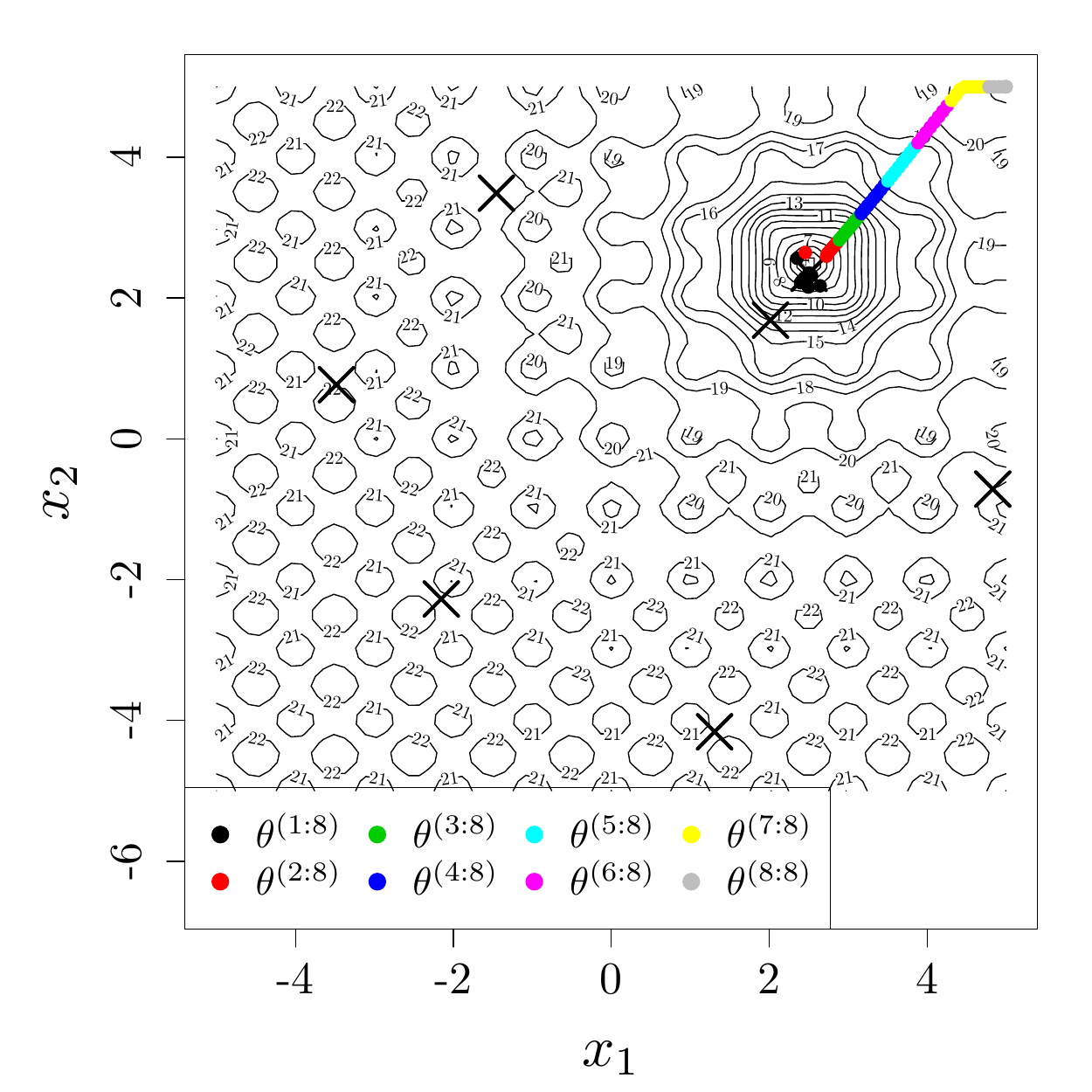}
\caption{First (top row) and second (bottom row) iteration of the toy greedy algorithm \ref{EGO_greedy} on the Sphere (left) and the Ackley functions(right). The initial DoE is different from the one shown in Figure \ref{x_theta_2D}. For more information see the caption of Figure \ref{x_theta_2D}.} 
\label{x_theta_2D_DoE2}
\end{figure}

\begin{figure}[H] 
\centering
\includegraphics[width=0.32\textwidth]{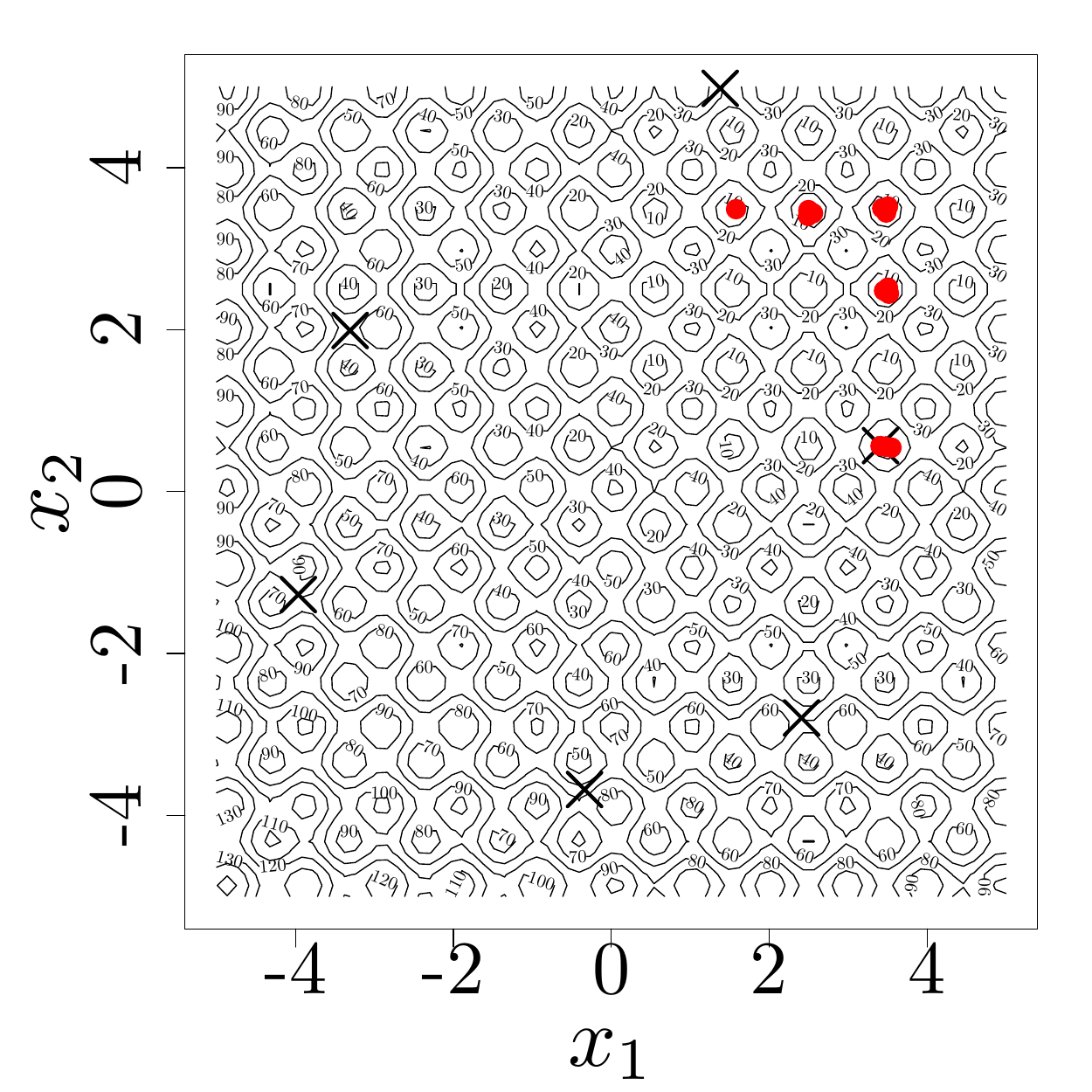}
\includegraphics[width=0.32\textwidth]{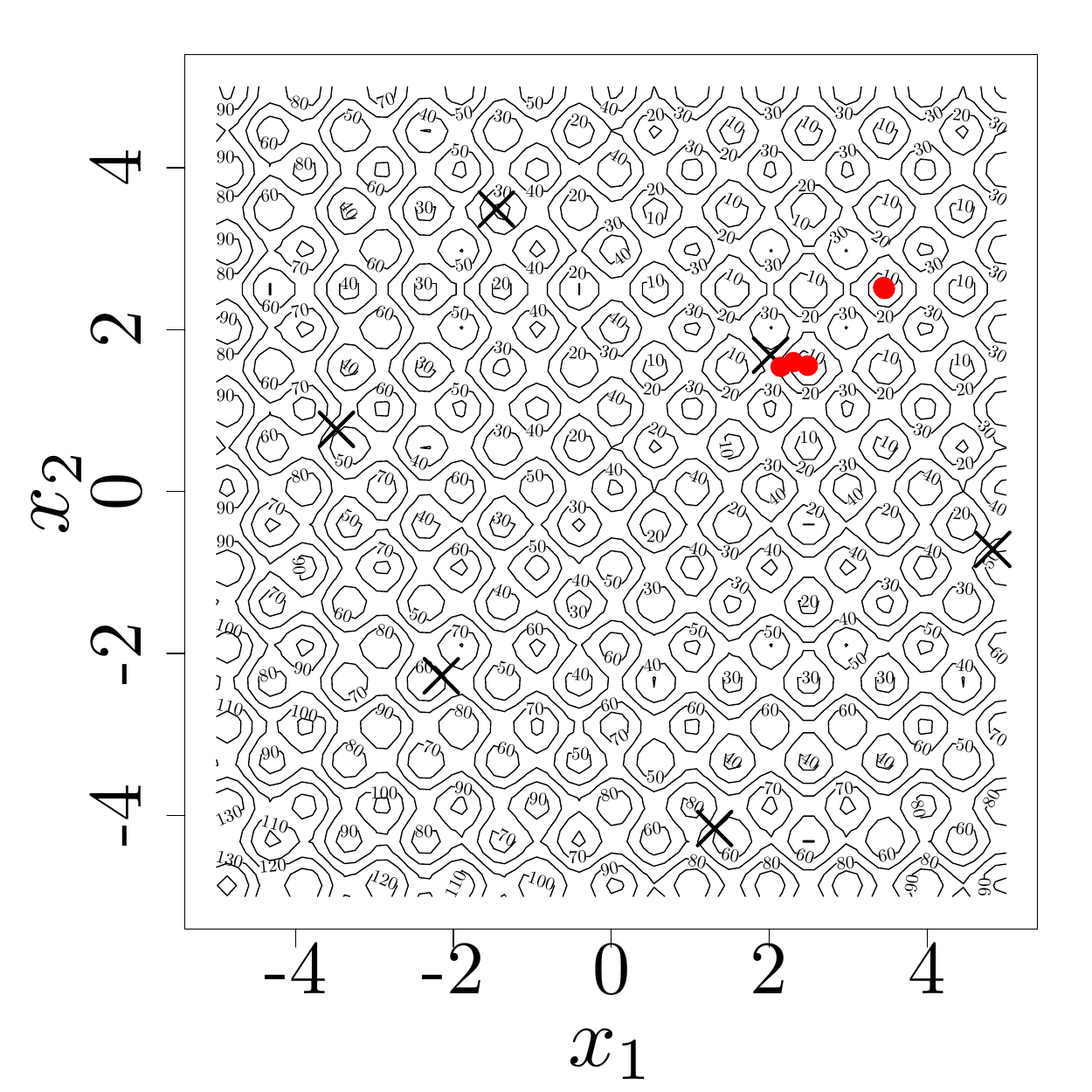}
\includegraphics[width=0.32\textwidth]{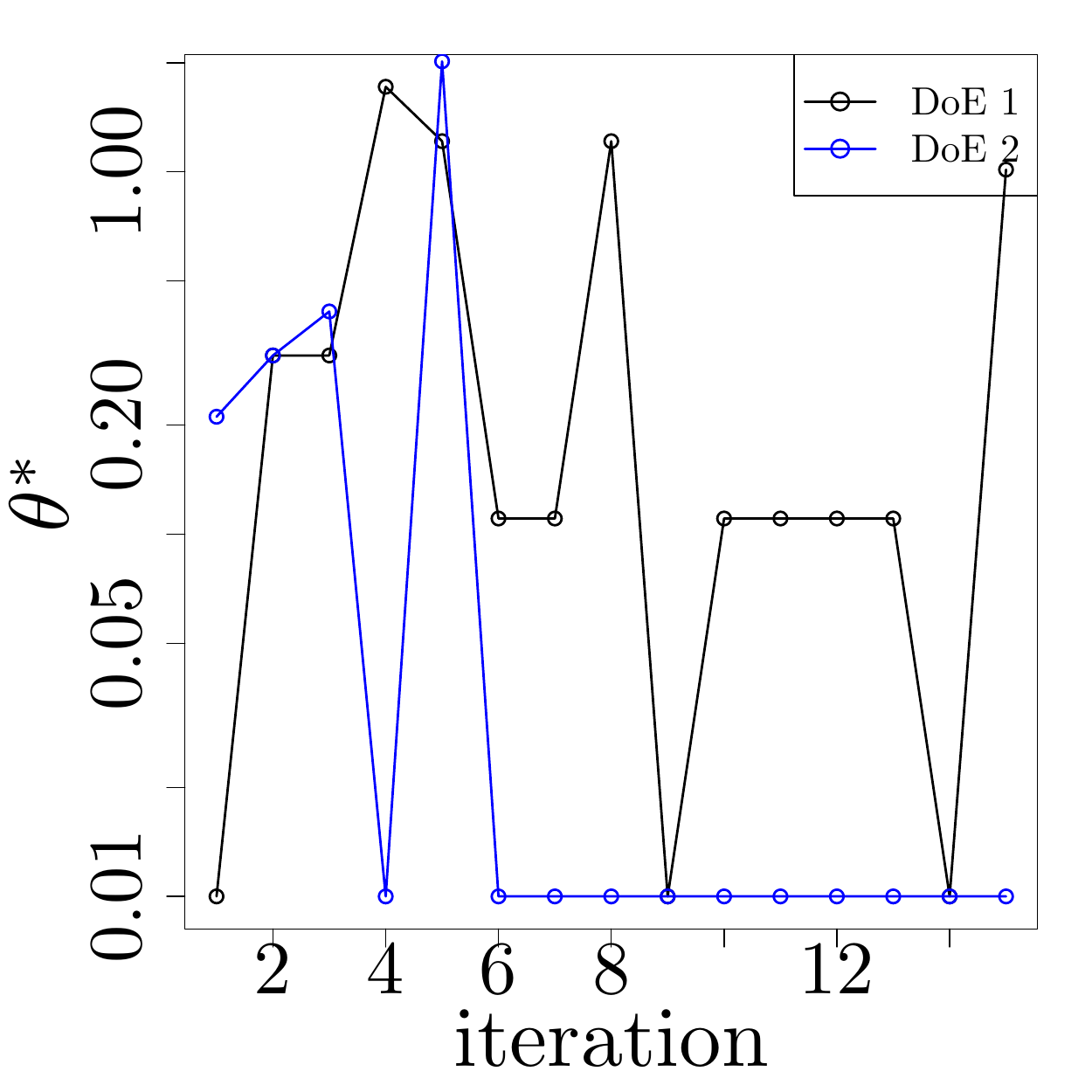}
\caption{DoEs created by the toy greedy algorithm \ref{EGO_greedy} after 15 iterations on the Rastrigin function with two DoEs (left and middle). Right: plot of ``best'' length-scale, $\theta^{*}$. The global minimum is located at $(2.5, 2.5)$.} 
\label{15_iter_rastrigin}
\end{figure}

In order to investigate the effect of initial DoE on the algorithm performance, the above experiments are repeated with another initial DoE.
Figure \ref{x_theta_2D_DoE2} shows the results which are similar to the previous experiments. For example, the length-scales tend to be small especially in the case of highly multimodal Ackley function. 
The algorithm's behavior, typical of small $\theta$'s (as explained in details in \cite{mohammadi2015}) is greedy,
that of a local search algorithm: local convergences can be seen in Figure \ref{15_iter_rastrigin} where the function to be optimized is Rastrigin with several local minima.
%%%%%%%%%%%%%%%%%%%%%%%%%%%%%%%%%%%%%%%%%%%%%%%%%%%%%%%%%%%%%%%%%%%%%%%%%%%%%%%%%%%%%%%%%%%%%%
%%%%%%%%%%%%%%%%%%%%%%%%%%%%%%%%%%%%%%%%%%%%%%%%%%%%%%%%%%%%%%%%%%%%%%%%%%%%%%%%%%%%%%%%%%%%%%
%%%%%%%%%%%%%%%%%%%%%%%%%%%%%%%%%%%%%%%%%%%%%%%%%%%%%%%%%%%%%%%%%%%%%%%%%%%%%%%%%%%%%%%%%%%%%%
\begin{figure}[H] 
\centering
\includegraphics[width=0.32\textwidth]{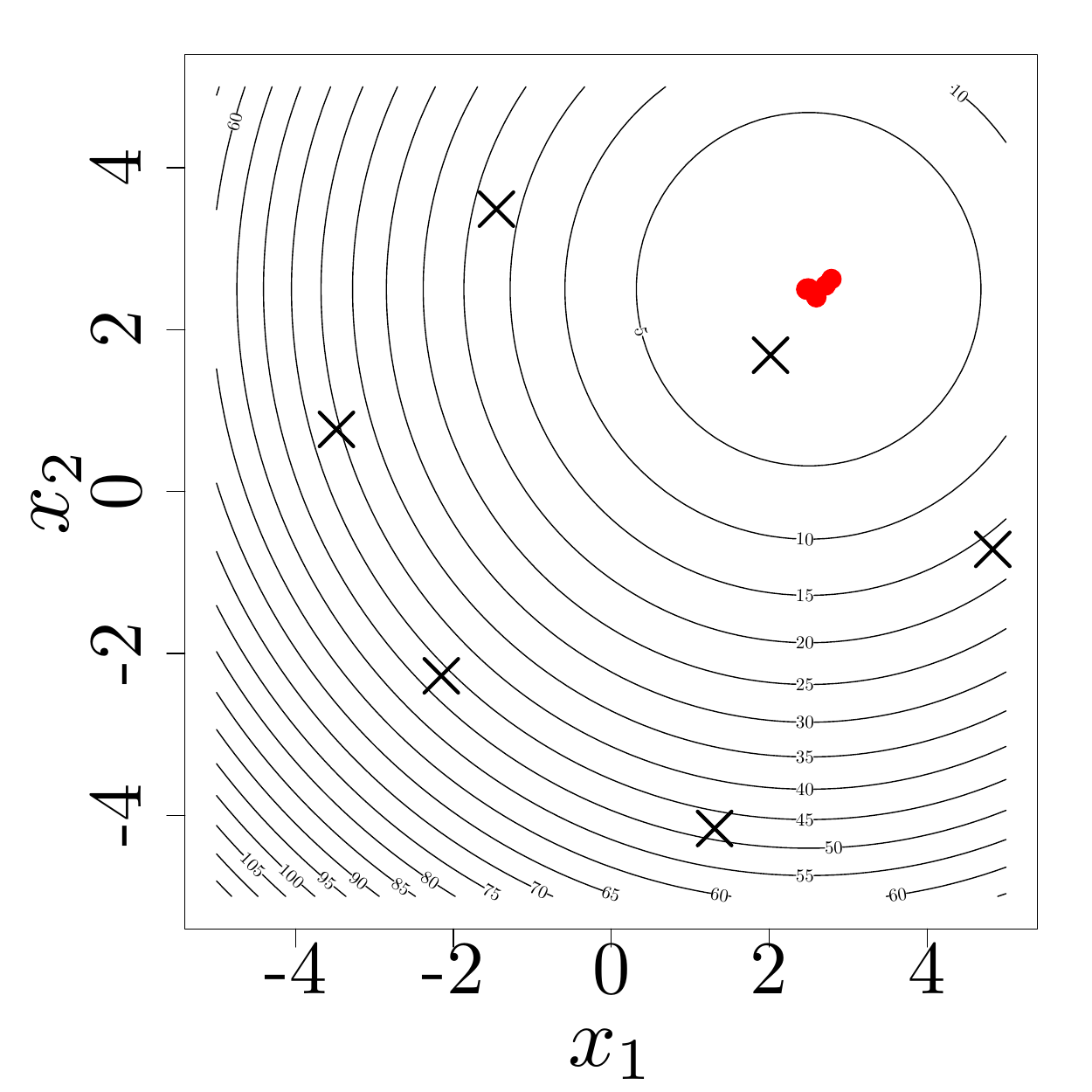}
\includegraphics[width=0.32\textwidth]{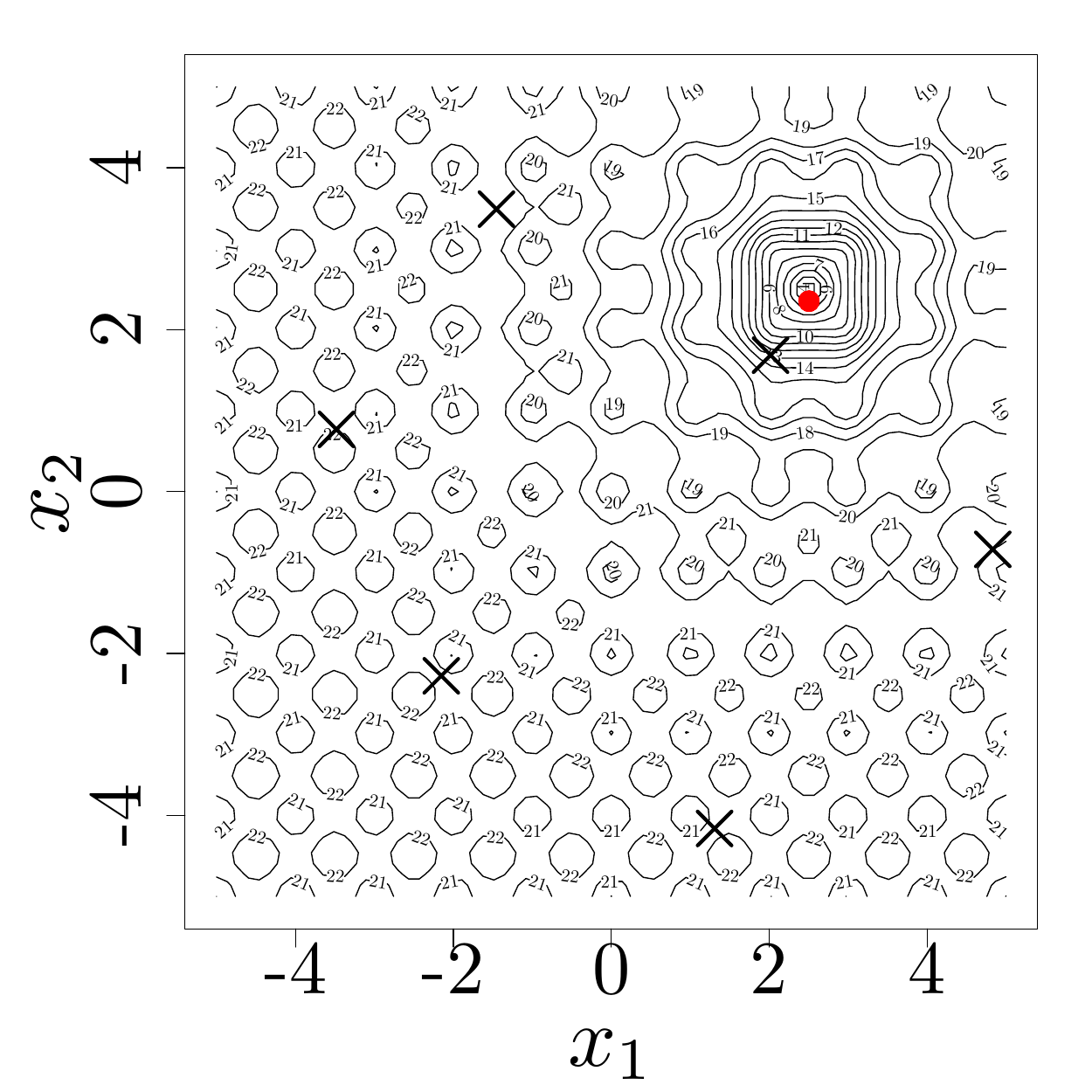}
\includegraphics[width=0.32\textwidth]{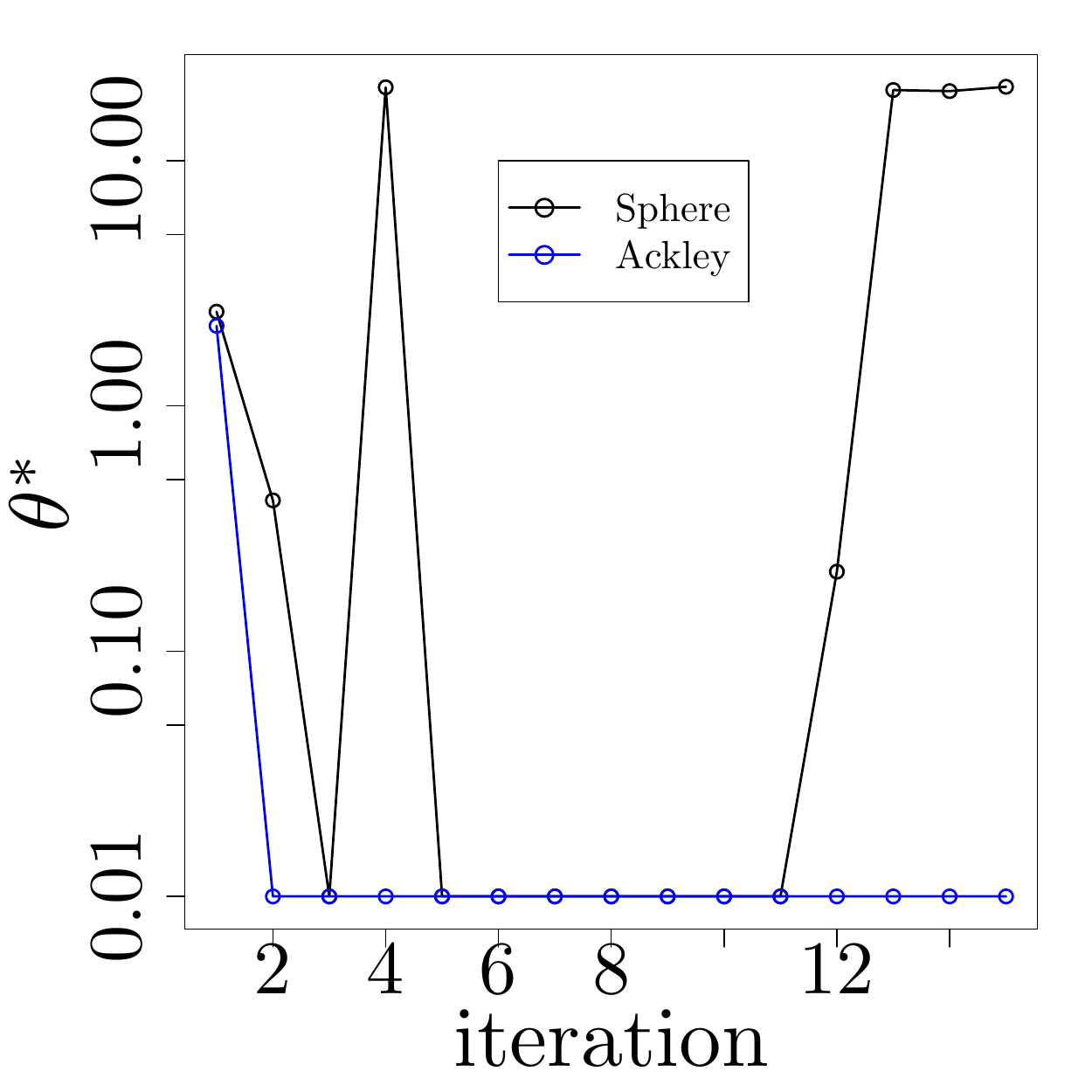}
\caption{DoEs created by the toy greedy algorithm \ref{EGO_greedy} after 15 iterations on the Sphere (left) and the Ackley (middle) functions. Right: plot of ``best'' length-scale, $\theta^{*}$. The initial DoE is different from the one shown in Figure \ref{15_iter}.} 
\label{15_iter_DoE2}
\end{figure}
%%%%%%%%%%%%%%%%%%%%%%%%%%%%%%%%%%%%%%%%%%%%%%%%%%%%%%%%%%%%%%%%%%%%%%%%%%%%%%%%%%%%%%%%%%%%%%

%%%%%%%%%%%%%%%%%%%%%%%%%%%%%%%%%%%%%%%%%%%%%%%%%%%%%%%%%%%%%%%%%%%%%%%%%%%%%%%%%%%%%%%%%%%%%%
%%%%%%%%%%%%%%%%%%%%%%%%%%%%%%%%%%%%%%%%%%%%%%%%%%%%%%%%%%%%%%%%%%%%%%%%%%%%%%%%%%%%%%%%%%%%%%

%=================================================================================================================================================================
\section{An EGO algorithm with a small ensemble of kriging models}
%=================================================================================================================================================================
\subsection{Description of the algorithm}
EGO is used for the optimization of computationally intensive functions. So, it is practically impossible to calculate $f \left(\textbf{x}^{n+1} (\theta| \textbf{X}) \right)$ for many length-scales in order to obtain $\theta^*$. 
Herein, we propose an approach that works with a limited number of kriging models. The ensemble 
 of kriging models is structured by the length-scales. 
The pseudo-code is given below (Algorithm \ref{EGO_ensemble_algorithm}) followed by a detailed explanation of the approach.        
\begin{algorithm}
\caption{EGO based on a small ensemble of kriging models}
\label{EGO_ensemble_algorithm}
\begin{algorithmic}
\STATE Create an initial design: $\textbf{X} = \left[ \textbf{x}^1, \dots , \textbf{x}^n \right]^\top$.
\STATE Evaluate function at $\textbf{X}$ and set $\textbf{y} = f(\textbf{X})$.
\STATE Set the maximum number of evaluations, $t_{\text{max}}$.
\FOR{t $~\leftarrow ~ $n+1 \TO $t_{\text{max}}$} 
\STATE Define a neighborhood of radius $R^{(t)}$ around the current sample points. 
\STATE Set $\textbf{X}^{(n+1)} = \emptyset$ and $\textbf{X}^{\text{sel}} = \emptyset$.
\STATE Generate $q$ length-scales, $\theta_1, \dots , \theta_q$.
\FOR{i $~\leftarrow ~ $1 \TO q}
\STATE $\textbf{x}^{n+1} ~\leftarrow ~ \arg\max_{\textbf{x} \in \mathcal S}~EI(\textbf{x}; \theta_i)$.
\STATE $\textbf{X}^{(n+1)} ~\leftarrow ~ \textbf{X}^{(n+1)} \cup \textbf{x}^{n+1}$.
\IF {$\textbf{x}^{n+1}$ is not inside the defined neighborhoods}
\STATE $\textbf{X}^{\text{sel}}~\leftarrow ~\textbf{X}^{\text{sel}} \cup \textbf{x}^{n+1}$.
\ENDIF
\ENDFOR
\IF {$\textbf{X}^{\text{sel}}=\emptyset$}
\STATE $\textbf{X}^{\text{sel}}~\leftarrow ~ \arg\max \left(\min_{\textbf{x} \in \textbf{X}^{(n+1)}} dist(\textbf{x}, \textbf{X})\right)$
\ENDIF
\STATE Evaluate function at $\textbf{X}^{\text{sel}}$ and set $\textbf{y}^{\text{sel}}=f(\textbf{X}^{\text{sel}})$.
\STATE Select $\theta^*$, for which $f(\arg\max_{\textbf{x} \in \mathcal S}~EI(\textbf{x}; \theta^*))=\min(\textbf{y}^{\text{sel}})$.
\STATE Generate two length-scales close to $\theta^*$. This yields two new infill samples by $EI$ maximization, $\textbf{X}^{\text{new}}=\left[\textbf{x}^{\text{new1}}, \textbf{x}^{\text{new2}} \right]^\top$. 
\STATE Evaluate function at $\textbf{X}^{\text{new}}$ and set $\textbf{y}^{\text{new}} = f(\textbf{X}^{\text{new}})$.
\STATE Update the DoE: $\textbf{X}~\leftarrow ~ \textbf{X} \cup \textbf{X}^{\text{sel}} \cup \textbf{X}^{\text{new}}$, $\textbf{y}~\leftarrow ~ \textbf{y} \cup \textbf{y}^{\text{sel}} \cup \textbf{y}^{\text{new}}$.
\ENDFOR
\end{algorithmic}
\end{algorithm}

Let $(\textbf{X}, \textbf{y})$ be the initial design of experiments. The covariance function we use here is the isotropic Mat\'ern 5/2 kernel \cite{GPML}. Thus, there exists only one length-scale to be tuned. The first reason for using an isotropic kernel is simplicity and clarity in the analysis. By taking isotropic functions and kernels, a difficult aspect of the algorithm (anisotropy, which is related to variables sensitivity) is neutralized to focus on other (also quite complex) phenomena. By taking isotropic kernels, the results of the numerical experiments are more stable. The second reason is that isotropic kernels have been found to perform well for EGO in high-dimension in the context of expensive-to-evaluate functions \cite{hutter2013}. 

At each iteration, five length-scales are generated. They are sampled on a basis 10 logarithmic scale from $[-2, 1]$ based on a Latin Hypercube Sampling (LHS) plan (that is $\theta$ ranges from $10^{-2}$ to $10^{1}$). Then, they are sorted and scaled back, $\theta_i=10^{\log\theta_i}~,~ 1\leq i \leq 5~;~\theta_1<\theta_2<\dots<\theta_5$. 
Corresponding to each length-scale $\theta_i$, a kriging model is created which gives a new infill sample: $\textbf{x}^{n+1}(\theta_i|\textbf{X})=\arg\max_{\textbf{x} \in \mathcal S}~EI(\textbf{x}; \theta_i)$. 

\begin{figure}[H] 
\centering
\includegraphics[width=0.5\textwidth]{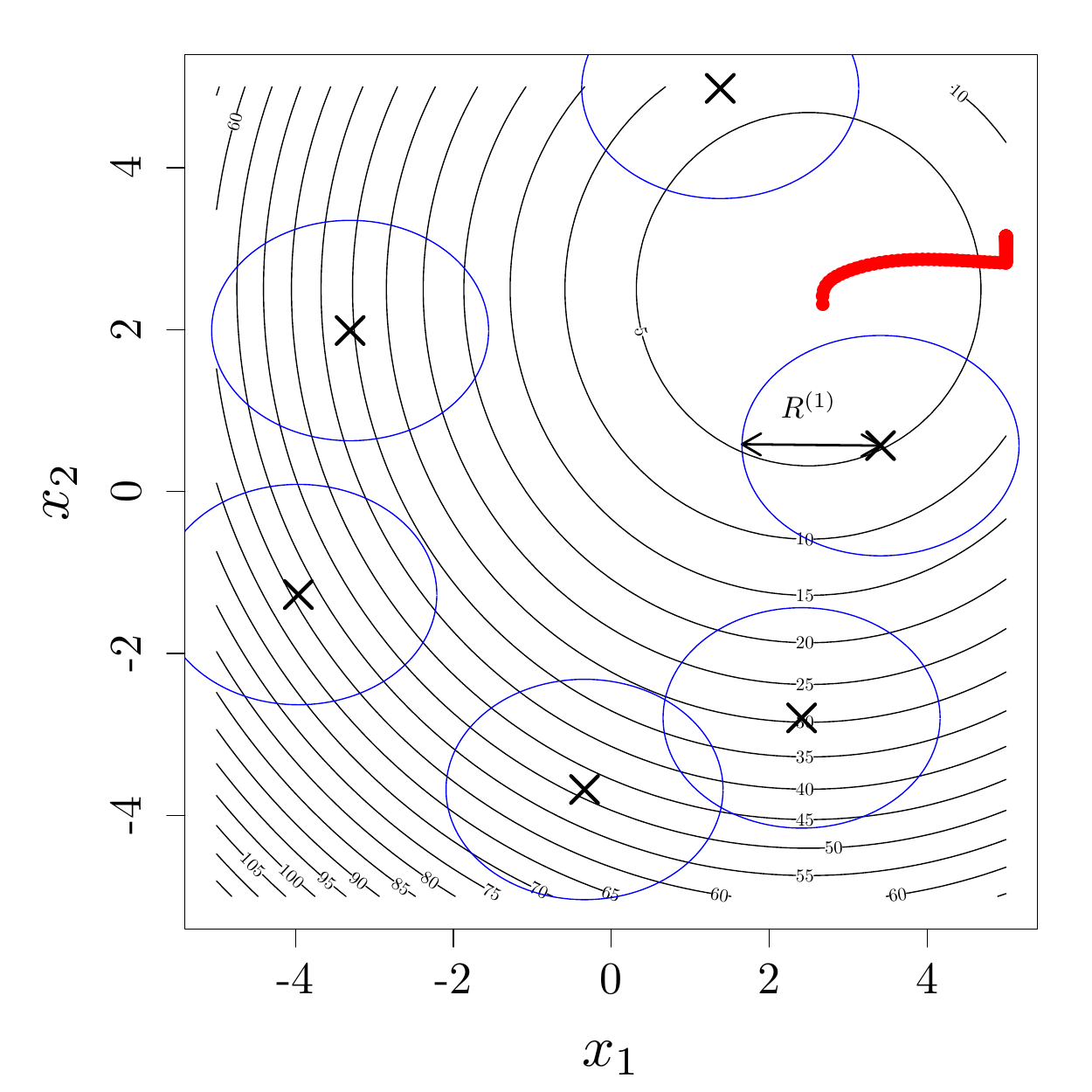}
\caption{DoE and neighborhoods as balls around the design points (blue circles). 
The infill samples occurring inside any neighborhood are not considered by the optimizer.} 
\label{EnsembleEGO_1st_iter}
\end{figure}

In the next step, the $\textbf{x}^{n+1}(\theta_i|\textbf{X})~,~1 \leq i \leq 5,$ that are not close to the design points are selected and the function is evaluated there. 
The notion of closeness is expressed by defining a neighborhood of radius $R^{(t)}$ around design points, see Figure \ref{EnsembleEGO_1st_iter}. 
It is important to prevent the points from converging around early good performers, otherwise such greedy algorithm where decisions are taken solely on the account of objective function values 
would not be sufficiently explorative for global optimization.
Further explanations about the neighborhood definition are provided in the next paragraph. 
The eligible $\textbf{x}^{n+1}(\theta_i|\textbf{X})~,~1 \leq i \leq 5,$ are selected and stored in the matrix $\textbf{X}^{\text{sel}}$. $\textbf{y}^{\text{sel}}$ contains the function values at $\textbf{X}^{\text{sel}}$.

%%%%%%%%%%%%%%%%%%%%%%%%%%%%%%%%%%%%%%%%%%%%%%%%%%%%%%%%%%%%%%%%%%%%%%%%%%%%%%%%%%%%%%%%%%%%%%
%%%%%%%%%%%%%%%%%%%%%%%%%%%%%%%%%%%%%%%%%%%%%%%%%%%%%%%%%%%%%%%%%%%%%%%%%%%%%%%%%%%%%%%%%%%%%%

The neighborhood defined around every design point is a ball with radius $R^{(t)}$ where the index $t$ is the iteration. 
As the optimization progresses, the radius shrinks according to the following linear scheme:
\begin{eqnarray}
R^{(t)}=
\left\{
	\begin{array}{ll}
		 R^{(1)} - \frac{R^{(1)}}{t_{\text{threshold}}}\times (t-1) & \mbox{if} ~~ t ~\leq ~t_{\text{threshold}} \\
		0 & \mbox{otherwise} ,
	\end{array}
\right. 
\label{Radius_neighbor}
\end{eqnarray}
in which $t_{\text{threshold}}$ is $70\%$ of total number of iterations, $t_{\text{max}}$. 
The initial radius $R^{(1)}$, is half of the distance between the best initial DoE (based on its $f$ value) and the closest design point to it. 
Again, defining such neighborhoods prevents the algorithm from focusing around good points too early.

Now, among the five generated length-scales, the best one is selected and is denoted by $\theta^*$. 
Recall that the best length-scale is the one that yields $f\left(\textbf{x}^{n+1}(\theta_i|\textbf{X})\right)=\min(\textbf{y}^{\text{sel}})$.
Then, two length-scales, $\theta_-^*$ and $\theta_+^*$, close to $\theta^*$ are generated. They are defined as:
\begin{itemize}
\item If $\theta^*=\theta_i, 2\leq i \leq 4$,  $\theta_-^*=\theta^* - \frac{1}{3}(\theta^* - \theta_{i-1})$ and $\theta_+^*=\theta^* + \frac{1}{3}(\theta_{i+1} - \theta^*)$. 
\item If $\theta^*=\theta_1$,  $\theta_-^*=0.01$ and $\theta_+^*=\theta^* + \frac{1}{3}(\theta_{2} - \theta^*)$.
\item If $\theta^*=\theta_5$,  $\theta_-^*=\theta^* - \frac{1}{3}(\theta^* - \theta_{4})$ and $\theta_+^*=10$.
\end{itemize}

The two new infill samples obtained with the kriging models with length-scales $\theta_-^*$ and $\theta_+^*$ are stored in the $\textbf{X}^{new}$ matrix, 
\begin{equation}
\textbf{X}^{\text{new}}= \left[\textbf{x}^{n+1}(\theta_-^*|\textbf{X}) ~,~ \textbf{x}^{n+1}(\theta_+^*|\textbf{X}) \right]^\top.
\end{equation}
Finally, the current DoE $(\textbf{X}$, $\textbf{y})$ is updated by adding $\textbf{X}^{\text{new}}$ and $\textbf{X}^{\text{sel}}$ to $\textbf{X}$ and $\textbf{y}^{\text{new}}$ and $\textbf{y}^{\text{sel}}$ to $\textbf{y}$. This procedure continuous until the budget is exhausted.

\subsection{Tests of the algorithm}
\begin{figure}[H] 
\centering
\includegraphics[width=0.49\textwidth]{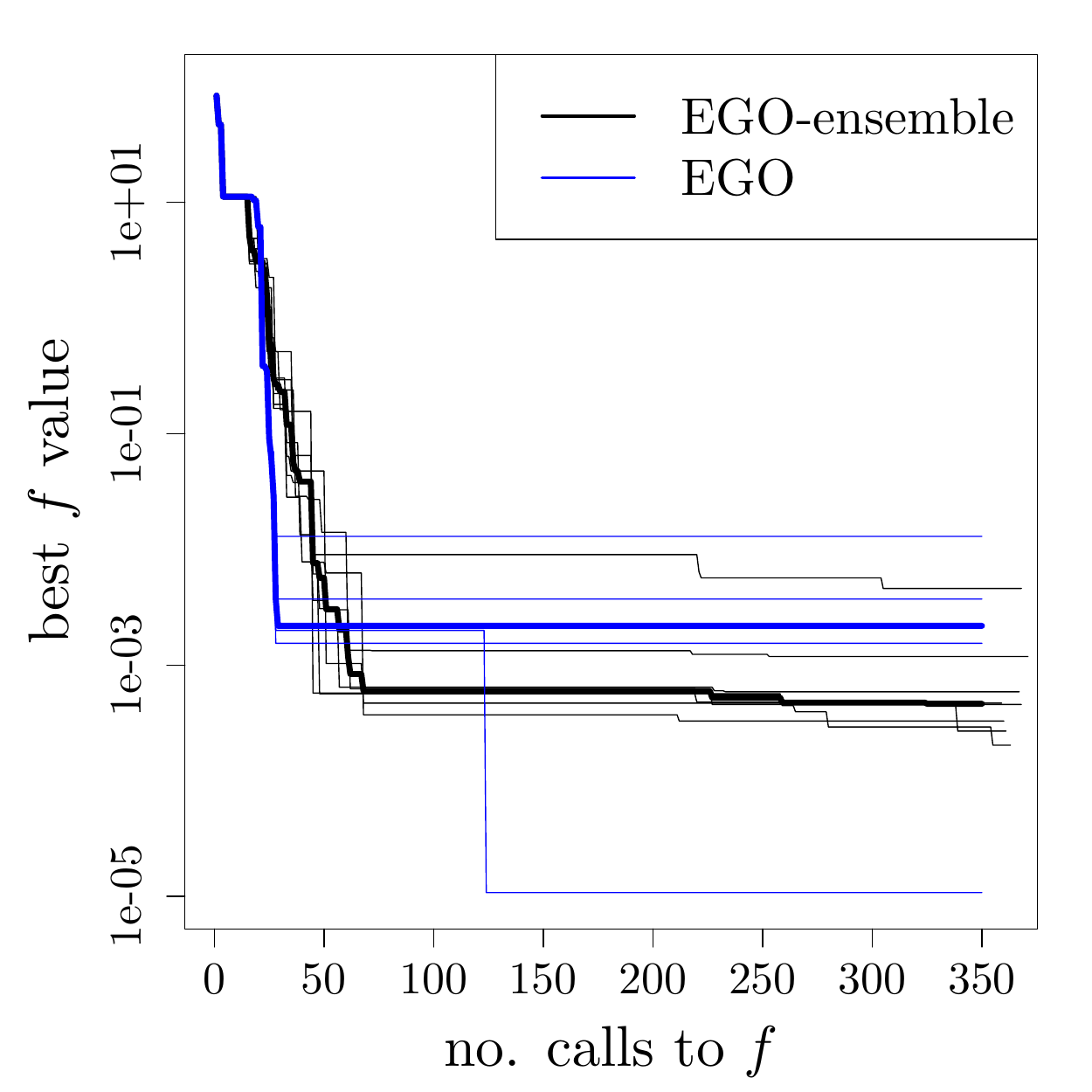}
\includegraphics[width=0.49\textwidth]{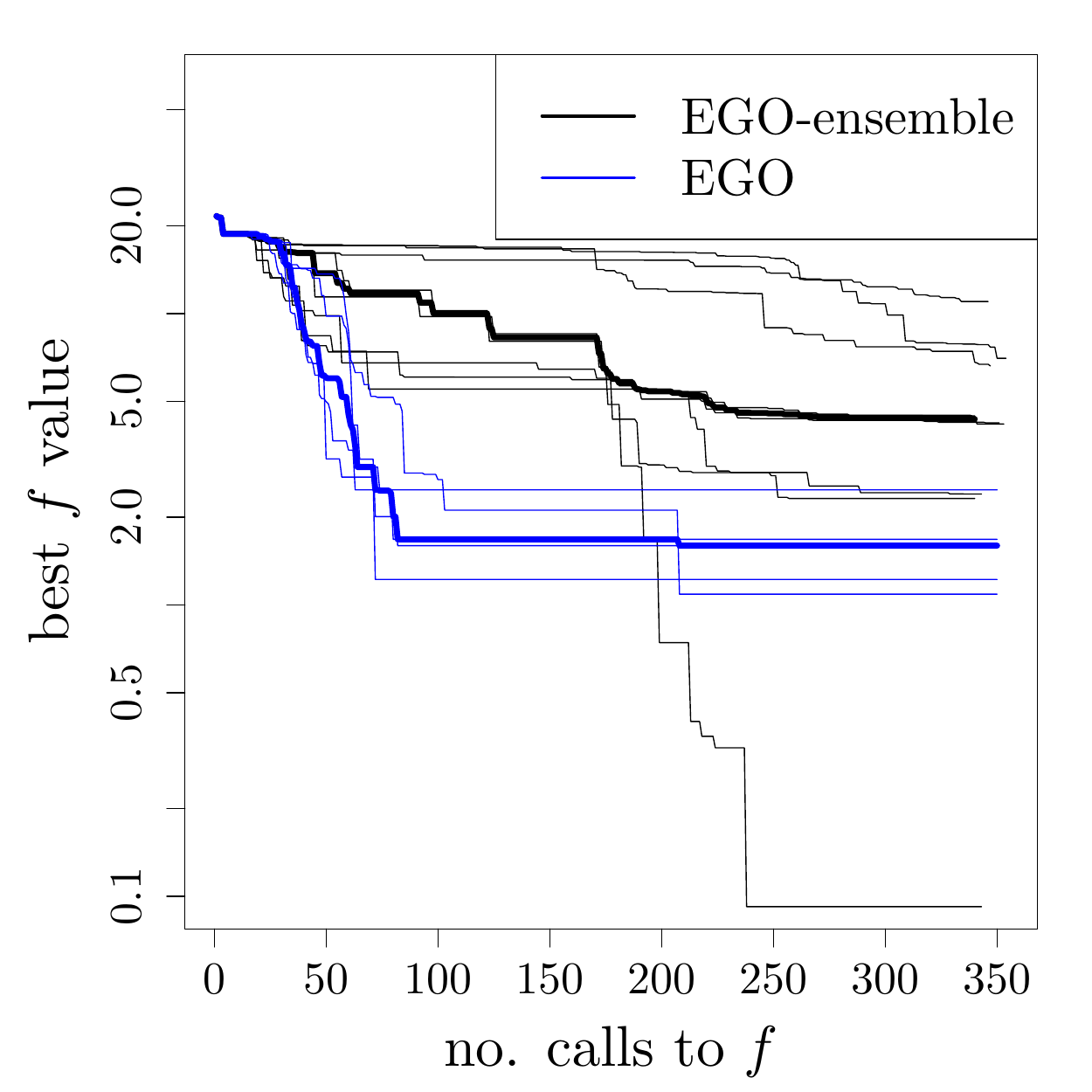}
\includegraphics[width=0.49\textwidth]{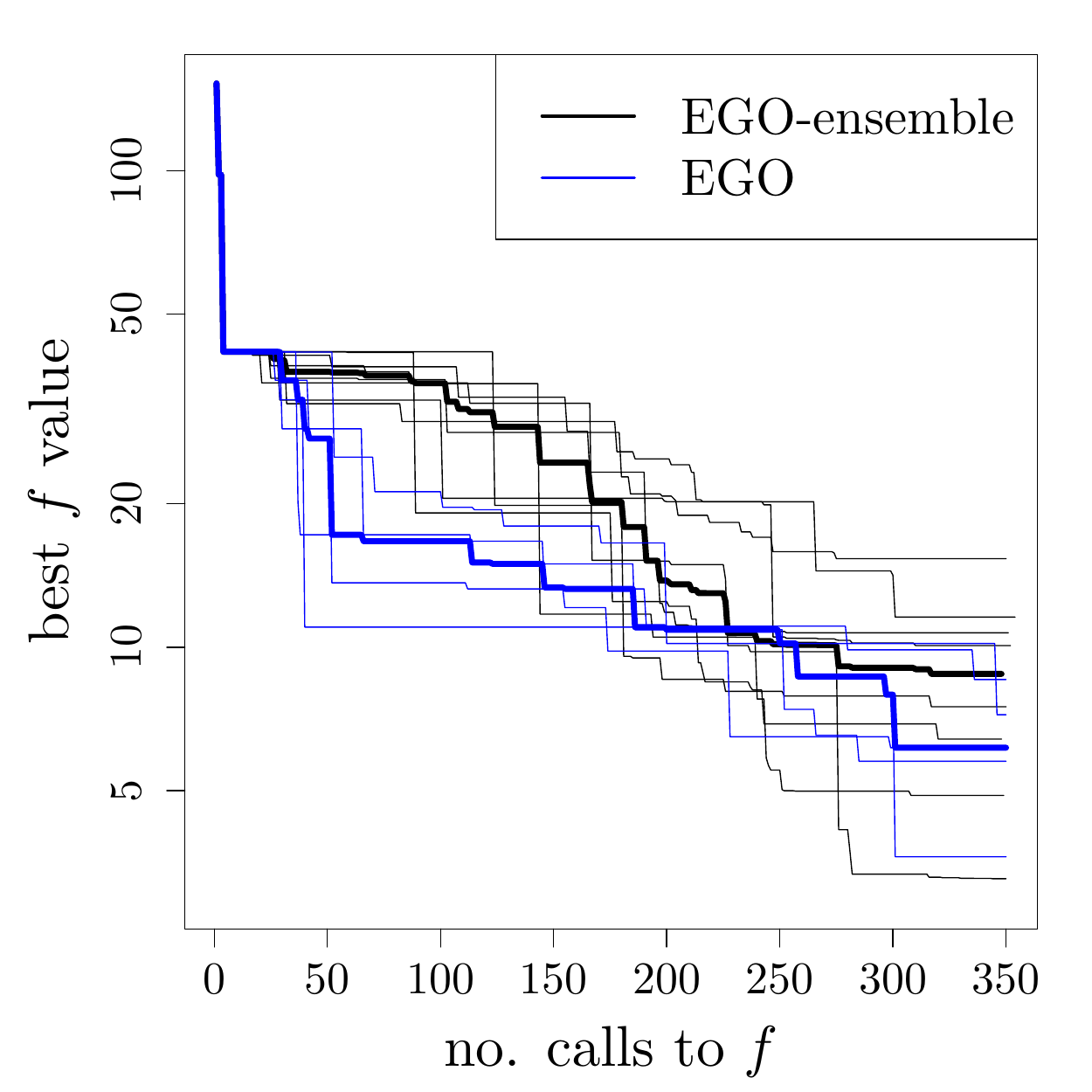}
\caption{Best objective function vs. number of calls of EGO with the ensemble of kriging models (thin black lines) and standard EGO (thin blue lines) on Sphere(top left), Ackley (top right) and Rastrigin (bottom) functions. The thick lines show the median of the runs.} 
\label{EnsembleEGO_EGO_dim5}
\end{figure}

The performance of this EGO method that is based on a small ensemble of kriging models (5+2 models) is tested on three isotropic functions, Sphere, Ackley and Rastrigin. 
The functions are defined in $\mathcal{S} = [-5, 5]^d$ where $d=5$. The total number of iterations is $15\times d$. Each optimization run is repeated eight times (thin black lines). 
Figure \ref{EnsembleEGO_EGO_dim5} shows the results and the performance of the standard EGO method (thin blue lines) which is repeated five times with a budget equals to $70\times d$. 
The plots show the best objective functions observed so far. 
The initial DoE is fixed for both algorithms and has a size equal to $3\times d$. 
The thick lines are the median of the runs.

The small ensemble version of EGO is slightly better on the sphere function because it benefits from 
its greedy choice of points that are never misleading. On Rastrigin and Ackley, the small ensemble EGO is slower early in the search, which might be due to the schedule of $R^{(t)}$. Later on, still on Rastrigin and Ackley, EGO with a small ensemble shows both the worst and best performances, therefore illustrating a tendency to get trapped in local optima. In terms of median performance, after 250 evaluations of the objective function (at the time when the neighborhood control ceases), the small ensemble EGO is equivalent to EGO on Rastrigin and worse on Ackley.

%%%%%%%%%%%%%%%%%%%%%%%%%%%%%%%%%%%%%%%%%%%%%%%%%%%%%%%%%%%%%%%%%%%%%%%%%%%%%%%%%%%%%%%%%%%%%%
%%%%%%%%%%%%%%%%%%%%%%%%%%%%%%%%%%%%%%%%%%%%%%%%%%%%%%%%%%%%%%%%%%%%%%%%%%%%%%%%%%%%%%%%%%%%%%
%=================================================================================================================================================================
\section{Conclusions}
We have investigated a variant of the EGO optimization algorithm where, instead of using at each 
iteration a kriging model learned through a statistical estimation procedure such as maximum likelihood, a small set of models with a fixed length-scale is employed. 
The motivations are threefolds. 
Firstly, it has been noticed in two-dimensions that the manifolds of the points that maximize expected improvement for various length-scales approach rapidly the global optimum. 
Secondly, ensemble methods have a lower computational complexity since the number of kriging covariance matrices inversions is limited to the number of elements in the ensemble, seven in the current work. On the contrary, maximum likelihood or cross-validation approaches require the inversion of the covariance matrix at each of their internal iteration.
Thirdly, ensemble methods may more easily lead to parallel versions of EGO as the maximization of expected improvement can be distributed on several computing nodes, one for each kriging model. 

Our first investigations have led to the following conclusions: tuning the length-scale to achieve an immediate improvement in the objective function may not be as efficient a strategy as two-dimensional plots of the manifold seem to indicate; the greediness of the method is a source of premature convergence to good performing points; optimal values of the length scale (in the sense of short term improvement) change a lot from one iteration to the next as the design of experiments evolves, rendering self-adaptive and Bayesian strategies not efficient for this purpose.

Nevertheless, we believe that the idea of searching in the space of length-scales as a proxy for searching in the space of optimization variables deserves further investigations because of its potential for tackling the curse of dimensionality. In particular, the schedule of the neighborhood radius, an iteration-smoothing learning procedure for the length-scales, and alternative strategies for making the ensemble of kriging should be studied. 

%=================================================================================================================================================================
\bibliography{biblio}
\bibliographystyle{plain}
%=================================================================================================================================================================
\end{document}